\documentclass[12pt]{article}

\usepackage{theorem}
\usepackage{amssymb}
\theoremheaderfont{\scshape}

\newtheorem{Proposition}{Proposition}
  \newtheorem{Remark}[Proposition]{Remark}
  \newtheorem{Corollary}[Proposition]{Corollary}
  \newtheorem{Lemma}[Proposition]{Lemma}
  
  \newtheorem{Theorem}[Proposition]{Theorem}

\def\Tm{\begin{Theorem}\label}
\def\Pp{\begin{Proposition}\label}
\def\Rm{\begin{Remark}\label}
\def\Lm{\begin{Lemma}\label}
\def\Co{\begin{Corollary}\label}

\def\eTm{\end{Theorem}}
\def\ePp{\end{Proposition}}
\def\eRm{\end{Remark}}
\def\eLm{\end{Lemma}}
\def\eCo{\end{Corollary}}

\def\Qp{{\bf Q}_\perp}

\def\FP{{\bf F}_\perp}
\def\dpp{{\bf D}_\perp}

\def\Rp{{\bf R}_\perp}

\makeatletter
\@addtoreset{equation}{section}
\makeatother

\def\Box{{\hfill\hbox{\enspace${\sqre}$}} \smallskip}
\def\sqr#1#2{{\vcenter{\vbox{\hrule height .#2pt
                             \hbox{\vrule width .#2pt height#1pt \kern#1pt
                                   \vrule width .#2pt}
                             \hrule height .#2pt}}}}
\def\sqre{\mathchoice\sqr54\sqr54\sqr{4.1}3\sqr{3.5}3}

\newcommand{\bT}[1]{\begin{theorem}\label{#1}}
\newcommand{\be}[1]{\begin{equation}\label{#1}}
\newcommand{\ba}[1]{\begin{eqnarray}\label{#1}}
\newcommand{\ee}{\end{equation}}
\newcommand{\ea}{\end{eqnarray}}
\newcommand{\bl}[1]{\begin{lemma}\label{#1}}
\newcommand{\bp}[1]{\begin{proposition}\label{#1}}
\newcommand{\br}[1]{\begin{remark}\label{#1}}

\def\n{{\cal N}}

\def\RR{\mathbb{R}}
\def\DD{\mathbb{D}}
\def\CC{\mathbb{C}}
\def\NN{\mathbb{N}}

\def\z{\noindent}
\def\bff{{\bf f}}
\def\bffz{{\bf F}_0}
\def\bfd{\mathbf{D}}

\def\bfh{{\bf h}}

\def\bfH{{\bf H}}\def\bfY{{\bf Y}}

\def\bfl{{\bf l}}

\def\bog{{\bf g}}
\def\bogl{{\bf g}^{(\bfl)}}
\def\bfgl{{\bf G_{l}}}
\def\calm{{\cal M}}

\def\calv{{\DD_\epsilon}}
\def\caln{{\cal N}}
\def\calnb{{{\cal N}}}
\def\lloc{{L}^1_{loc}}
\def\lray{{L_{ray}}}
\def\lb{L^1_b}
\def\lone{{L^1}}

\def\ga{\hat{\Gamma}}
\def\gc{\hat{\Gamma}_c}
\def\hatc{\hat C}
\def\bc{\hat{B}_c}

\def\hb{\hat{B}}

\def\bfa{{\bf a}}
\def\bfB{{\bf B}}
\def\bfA{{\bf A}}
\def\bfQ{{\bf Q}}
\def\zpp{(z\oplus\gc(z))}

\def\betapp{(\beta\oplus\bc)}

\def\bfy{{\bf y}}
\def\hatby{\tilde{\bfy}}

\def\bfR{{\bf R}}
\def\bfz{{\bf z}}
\def\lap{{\cal L}}
\def\lapi{{\cal L}^{-1}}
\def\bfZ{{\bf Z}}
\def\bfk{{\bf k}}

\def\bfm{{\bf m}}
\def\bfdl{{\bf d_l}}
\def\bfdm{{\bf d_m}}
\def\lap{ }
\def\lap{{\cal L}}
\def\lapi{{\cal L}^{-1}}
\def\bor{{\cal B}}
\def\cald{{\cal D}}

\def\bfZ{{\bf Z}}
\def\bfk{{\bf k}}

\def\bfm{{\bf m}}
\def\bfdl{{\bf d_l}}
\def\bfdm{{\bf d_m}}
\def\heav{{\cal H}}

\title { Exponential asymptotics, trans-series 
and generalized Borel summation for analytic nonlinear rank one systems
of
ODE's
}

\author{Ovidiu Costin \thanks {Mathematics Department, Hill Center
Rutgers University, New Brunswick, NJ 08903;
e-mail: costin@maxwell.rutgers.edu}}
\date{ }

\begin{document}
%\mbox{ }
%\eject\addtocounter{page}{-1}

\include{psfig}
\maketitle
\begin{abstract}
For analytic nonlinear systems of ordinary differential equations, under
some non-degeneracy and integrability conditions we prove that 
the formal exponential series
solutions (trans-series) at an
irregular singularity  of
rank one are Borel summable (in a sense similar to that
of Ecalle).  The functions obtained by re-summation of the trans-series
are precisely the solutions of the differential equation that decay in a
specified sector in the complex plane.

We
find the dependence of the correspondence
between the solutions of the differential equation and trans-series
as the ray in the complex plane changes (local Stokes phenomenon).

We study, in addition, the general solution in $\lloc$ of the convolution
equations corresponding, by inverse Laplace transform,
to the given system of ODE's, and its analytic properties.

 Simple analytic identities
lead to ``resurgence'' relations and to an averaging formula
having, in addition  to the properties of the medianization
of Ecalle, the property of preserving exponential growth
at infinity.

\end{abstract}

\begin{section} {Introduction and main results}

We consider an $n$-dimensional, rank one, level-one
vector differential equation
in a neighborhood of  an irregular singularity, say $x=\infty$.
We assume that the Stokes lines are simple.  In normalized
form (see \cite{Wasow},
\cite{To1} ), such an equation can be written in the form

\be{eqor}
{\bf y}'={\bf f}_0(x)-\hat\Lambda {\bf y}-
\frac{1}{x}\hat B {\bf y}+{\bf g}(x,{\bf y}),\ \ \bfy\in\CC^n,\
\ee

(The reason to separate out the second and third term on the
r.h.s. of (\ref{eqor}) is that they play a special role in the
asymptotic behavior of the solutions).

The functions $\xi\mapsto{\bf f}_0(\xi^{-1})$
and $(\xi,\bfy)\mapsto\bog(\xi^{-1},\bfy)$ are taken to be
analytic for small arguments. The normalization can be chosen
so that $\mathbf{f}_0(x)=O(x^{-2})$ for large $x$, and,
by construction, we have $\mathbf{g}(x,\bfy)=O(|\bfy|^2,x^{-2}\bfy)$.

$\hat\Lambda$ and $\hat B$ are $n\times n$ matrices
with constant coefficients.
 We assume that $\hat\Lambda$ is invertible and that
the (``non-resonance'') condition
$\arg\lambda_j\ne\arg\lambda_i,$ for $\ j\ne i$, $\lambda\in$\,spec\,$\hat
\Lambda$, is satisfied.

By a change of variables we can then arrange that $\hat\Lambda$ is diagonal,
$\hat\Lambda={\mbox{diag}}\{\lambda_i\}$ with 
$\arg \lambda_j>\arg\lambda_i$ for $j>i$ and make
$\lambda_1=1$. The matrix $\hat B$ can be diagonalized at
the same time \cite{To1}.

To simplify the analysis we  assume further
that  $\Re(\beta)>0$ where $\beta=\hat B_{1,1}$. 
Through normalization we make

\be{defbeta}
\Re(\beta)\in(0,1]\ee

We are interested  in the study of  the
solutions of (\ref{eqor}) that are decaying for large $x$, in 
one of the half-planes $\Re (xe^{-i\phi})>0$ with
$\phi\in(\arg\lambda_n-2\pi,
\arg\lambda_2)$.
 These
solutions have the same asymptotic behavior
at large $x$, described by a (typically divergent) power series

\be{asymbh}
\bfy(x)\sim\hatby_0=\sum_{k=2}^{\infty}\frac{\tilde\bfy_{0,k}}{x^k}
\ \ (|x|\rightarrow\infty;\ \Re\left(\,xe^{-i\phi}\right)>const>0)
\ee

\z For instance, all the solutions of the equation $y'+y=x^{-1}$
have the property $y(x)\sim\sum_{k=0}^{\infty}k!x^{-k-1}$ as $x\rightarrow\infty$.
If $\phi\ne 0$ there is only one solution of (\ref{eqor})
satisfying (\ref{asymbh}).
A much more  interesting case is when we take $\phi=0$.
Then, 
as it is known (and  will also follow from the present paper) 
there is a one dimensional manifold $M^+$ of solutions of
(\ref{eqor}) such that (\ref{asymbh}) holds.
The manifold $\tilde M^+$
of all {\em formal} solutions which decay
in the half plane $\Re x>0$
\be{genforms}
\hatby=\hatby_0+\sum_{k=1}^{\infty}C^k e^{-kx}\hatby_k
\ee

\z  also has one free parameter, $C\in\CC$. In (\ref{genforms}), $\hatby_k$, $k\ge 0$, are formal
power series and 
$\hatby$ is an instance of a trans-series. In our example
$y'+y=x^{-1},\  \hatby=\sum_{k=0}^{\infty}k!x^{-k-1}+Ce^{-x}$.
See Section~\ref{cor-for} 
a  heuristic construction leading to trans-series
solutions and for references.

 The series $\hatby_k$
satisfy the system of differential equations

$$
\bfy_0'+\left(\hat\Lambda +
\frac{1}{x}\hat B\right) \bfy_0={\bf f}_0(x)
+{\bf g}(x,\bfy_0)\phantom{mmmmmmmmmmmmmmmmm}$$
\vskip -0.3cm
\be{systemform}
\bfy'_k+
\left(\hat\Lambda+\frac{1}{x}\hat B-k-
{\mathbf{\partial}}{\bf g}(x,\bfy_0)\right)
\bfy_k=\sum_{|\bfl|> 1}\frac{\bogl(x,\bfy_0)}{\bfl!}
\sum_{\Sigma m=k}\prod_{i=1}^n\prod_{j=1}^{l_i}(\bfy_{m_{i,j}})_i
\ee

\z where $\bogl:={\partial^{(\bfl)}{\bf g}}/{\partial\bfy^\bfl}$,
${(\mathbf{\partial}{\bf g} )\bfy_k}:=\sum_{i=1}^n
(\bfy_k)_i({\partial {\bf g}}/{{\partial y_i}})$, and
$\sum_{\Sigma m=k}$ stands for the sum over all integers 
$m_{i,j}\ge 1$ with $1\le i\le n,1\le j\le l_i$
such that $\sum_{i=1}^n\sum_{j=1}^{l_i}m_{i,j}=k$. Because
$m_{i,j}\ge 1$, $\sum m_{i,j}=k$ (fixed) and card$\{m_{i,j}\}=|\bfl|$,
 the sums in (\ref{systemform})
contain only a {\em finite} number of terms.
We use the convention $\prod_{i\in\emptyset}\equiv 0$. 
The system (\ref{systemform}) is derived in Section ~\ref{cor-for}.

Starting with $k=1$
 the equations (\ref{systemform}) are linear.  Note that the inhomogeneous
term in these linear equations is zero for $k=1$, and 
for $k>1$ it involves only  $\bfy_n$ with
$n<k$.

 While some  connection between 
$\hatby_0$ and actual  solutions 
of (\ref{eqor}) is given by (\ref{asymbh}), the interpretation
of (\ref{genforms}) is less immediate, since generically
all the series involved are (factorially) divergent
and ``beyond all orders of each other''. 
The interest in trans-series is motivated partly
by their formal simplicity compared to the
vast class of differential equations that 
they ``solve'' and by the fact
that they can be algorithmically found,
once the equation is given. 
Finding the connection between formal expansions
and true solutions is the object of exponential asymptotics,
a field that has been growing constantly, especially
after the pioneering works of M. Berry,
 J. Ecalle
and M. Kruskal. 

The formalism of generalized Borel summation as well as the theory of
trans-series, in a very comprehensive setting, were introduced by Ecalle
\cite{Ecalle-book}, \cite{Ecalle}, \cite{Ecalle2}.

For the problem (\ref{eqor})---(\ref{asymbh}), we prove that there is a one-to-one natural correspondence
between actual solutions $\bfy$ and the trans-series $\hatby$
(\ref{genforms}).

We show that the general solution
of (\ref{eqor}), (\ref{asymbh}) is obtained by replacing each formal series in 
(\ref{genforms}) by its Borel sum which gives  a one-to-one
correspondence
between the formal solutions (trans-series)  and the actual solutions
of (\ref{eqor}), (\ref{asymbh}):

\be{eqrep1}
\hatby=\hatby_0+\sum_{k=1}^{\infty}C^k e^{-kx}\hatby_k\
\longleftrightarrow\ \lap_\phi{\cal B}_\phi\hatby_0+\sum_{k=1}^{\infty}C^ke^{-kx}
\lap_\phi{\cal B}_\phi\hatby_k=\bfy\ee

\z The Borel summation operator, $\lap\cal B$ will
be defined precisely. The function $\bfy\in M^+$ is
 convergently defined by (\ref{eqrep1}) for large
$x$. The left arrow in (\ref{eqrep1}) means that  $\lap_\phi{\cal
B}_\phi\hatby_k(x)\sim\hatby_k(x)$ for  $x\rightarrow\infty$.
The exact statement corresponding to (\ref{eqrep1}) is given
in Theorem~\ref{teo2}.

We study in detail the features of the representation
(\ref{eqrep1}) and the properties of the objects involved. The technique
that we use differs from that of \cite{Ecalle-book}, 
\cite{Ecalle}, \cite{Ecalle2} and leads to new results. In particular
we obtain for the Borel transform of the formal series solutions
of differential systems an 
averaging formula, having, as the medianization of Ecalle the quality
of preserving  reality and of commuting with convolution, but involving
a smaller number of analytic continuations and in addition 
satisfying the condition of at most exponential growth at infinity.

  For $m>1$, the inverse Laplace transform of $x^{-m}$ is
 $$\lapi x^{-m}
= p^{m-1}/\Gamma(m-1)={\cal B}\,x^{-m}$$

The \emph{Borel transform} ${\cal B}$
of a formal series 

\be{defyk}
\hatby=x^r\sum_{k=1}^{\infty}\tilde\bfy_k x^{-k},\  r\in (0,1)
\ee

\z is by definition the formal series gotten  
 by taking $\lapi$ term by term:

\be{defborel}
{\cal B}\,\hatby=\bfY:=p^{-r}\sum_{k=0}^{\infty} \frac{\tilde
\bfy_{k+1}}{\Gamma(k-r)}p^{k}
\ee

A priori $\bfY$ is still a formal series. If it has a nonzero
radius of convergence, then it generates an element of an analytic
function which we will  denote, all the same, by $\bfY$.

A formal series $\hatby$ is Borel summable in the classical sense
along a ray $\Phi$ (the
direction of which is given by the angle $\phi$) if the following
conditions are met:

1) The series $\bfY$ has a nonzero radius of convergence;

2) $\bfY$ can be analytically continued along the ray and

3) The analytic continuation $\bfY$ grows at most exponentially along
the ray and is therefore Laplace transformable along $\Phi$.

The Laplace transform along that ray of $\bfY$, $\lap_\phi\bfY$, is
well defined and gives the so
called \emph{Borel sum} of $\hatby$. We prove that the conditions 1 through
3 are met by $\bor \hatby_k$, $k\ge 0$, away from the Stokes rays, i.e.,
if $\phi\ne\arg\,\lambda_i, \ \lambda_i\in$spec$\,\hat\Lambda$.

Of all the formal solutions (\ref{genforms}), only the one with $C=0$
(formally) decays in a half-plane, if the half-plane {\em is not}
centered on the real axis. On the other hand, $\lap_\phi{\cal
B}\hatby_0$ turns out to be the only solution of (\ref{eqor}),
(\ref{asymbh}) which decays in the same half-plane centered on $\Phi$.
Borel summation associates uniquely a true solution to $\bfY_0$.

 The situation is more complicated and more interesting along Stokes rays
(we focus on one of them, $\Phi=\RR^+$).  Condition
2) above is violated and, generically, the functions
$\bfY_k$ have an array of branch points along $\RR^+$. If we reinterpret
2) and consider paths that avoid the singularities then first of all,
analytic continuation is (a priori) ambiguous. What is worse, the Laplace
transform of such analytic continuations  of $\bfY_0$ are,
typically,  \emph{not} solutions of (\ref{eqor})
 (see Section~\ref{noncomut}). However,
Laplace transforms of (a one-parameter family) of suitable 
weighted combinations of analytic continuations of $\bfY_0$
\emph{are}, as we will prove, solutions of (\ref{eqor}). If  we require
in addition that real series are Borel-summed to real-valued functions
then one
of weighted average of analytic continuations
appears as more natural (see also Theorem~\ref{Tecal} below).

\centerline{*}

To define the Borel transform along the Stokes line $\RR^+$ we construct
a suitable space of analytic functions.  Let $\phi_+=\arg\lambda_2$,
$\phi_-=2\pi-\arg\lambda_n$, and \be{defsect}
\mathcal{W}_1:=\{p:p\not\in\NN\cup\{0\}\mbox{ and }\
\arg\,p\in(-\phi_-,\phi_+)\} \ee (Fig. 1), a sector containing only the
eigenvalue $\lambda_1=1$ and punctured at all the integers (where the
functions ${\cal B}\hatby_k$ are typically singular; if $n=1$ the
condition on the argument is dropped). We construct over $\mathcal{W}_1$
a surface $\mathcal{R}_1$, consisting of homotopy classes of curves
starting at the origin, going only forward and crossing the real axis at
most once:

\ba{defpaths}
&&{\cal R_{\mathrm 1}}:=\Big\{\gamma:(0,1)\mapsto \mathcal{W}_1\ \  \mbox{s.t.\ }
\gamma(0_+)=0;\ \Re\,(\gamma(t))\ \mbox{increases
in $t$}\ \mbox{and}\cr
&&
\ 0=\Im(\gamma(t_1))=\Im(\gamma(t_2))
\Rightarrow t_1=t_2\Big\}\cr
&&
\end{eqnarray}

\z modulo homotopies. Let also

\be{defD}
{\cal D}:=\CC\backslash \cup_{i=1}^n\{\alpha\lambda_i:\alpha\ge 1\}
\ee

\z  be the complex plane without 
the rays originating at the eigenvalues $\lambda_i$ of 
$\hat\Lambda$.

\newcounter{cms}
\setlength{\unitlength}{1mm}
\begin{picture}(50,70)
\put(15,40){\makebox(0,0)[bl]{\large{O}}}\thicklines
\multiput(31,38)(10,0){3}{\addtocounter{cms}{1}\makebox(0,0)[b]
{$\mathsf{\scriptstyle{\arabic{cms}}}$}}
\multiput(61,39)(10,0){1}{\addtocounter{cms}{1}\makebox(0,0)[b]
{$\mathsf{\scriptstyle{\arabic{cms}}}$}}
\multiput(71,38)(10,0){4}{\addtocounter{cms}{1}\makebox(0,0)[b]
{$\mathsf{\scriptstyle{\arabic{cms}}}$}}
%\put(15,20){\circle{6}}
%\put(15,20){\circle*{2}}
\thinlines
\put(20,40){\line(4,1){80}}
\put(20,40){\line(3,-1){80}}
\thicklines
\multiput(30,40)(10,0){8}{$\scriptscriptstyle \diamond$}\
\put(20,40){\qbezier[15](0,0)(8,-2)(15,-4)\qbezier[15](15,-4)(23,-6)(30,-4)
\qbezier[15](30,-4)(38,-2)(45,0)\thicklines \qbezier[15](45,0)(53,2)(60,4)
\qbezier[15](60,4)(68,6)(75,4)\put(75,4){\large$(-^{4}+)$}}
\put(55,30){\LARGE${\mathcal{W}}_1$}
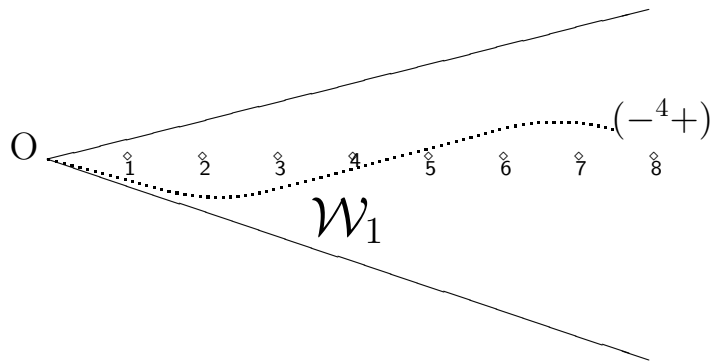
\put(0,05){{Fig 1.} \emph{The region $\mathcal{W}_1$. The dotted
line is one of the paths that generate $\mathcal{R}_1$.
}}
\end{picture}

Using notations similar to those of  Ecalle, we symbolize
the paths in $\mathcal{R}_1$ by  a sequence of signs
${\epsilon_1,..,\epsilon_j,..,\epsilon_n},\ \epsilon_j= +$ or $-$.
For example, $----+=-^4+$
will symbolize a path in  $\mathcal{R}_1$ that crosses
the real line from below through the interval $(4,5)$, and
then goes only through the upper half-plane (Fig.1);\ 
$''+''$ is a path confined to the upper half plane, etc. The analytic continuation of a function
$\bfY$ along the path $-^4+$ will be denoted $\bfY^{-^4+}$.

The result below gives a first characterization
of the analytic properties
of $\bor \hatby_k$. 
(In the following, we choose the determination of the logarithm which
is real for positive argument.)

\Pp{basicpp}

i) The function $\bfY_0:=\bor \hatby_0$ is analytic in $\cald$
and Laplace
transformable
along any direction in $\cald$. In a neighborhood of $p=1$

\be{caracty0p=1}
\bfY_0(p)=\left\{\begin{array}{l}
S_\beta(1-p)^{\beta-1}\bfA(p)+\bfB(p)\ \mbox{for $\beta\ne 1$}\cr
S_\beta\ln(1-p)\bfA(p)+\bfB(p)\ \mbox{for $\beta=1$}\end{array}
\right.
\ee

\z (see (\ref{defbeta})), where 
$\bfA$, $\bfB$ are ($\CC^n$-valued) analytic functions 
in a neighborhood of $p=1$.

\bigskip

ii) 
 The functions
$\bfY_k:=\bor \hatby_k,\ k=0,1,2,..$ are analytic in ${\cal R_{\mathrm
1}}$.

iii) For small $p$,

\be{analiticstructn}
\bfY_0(p)=p\bfA_0(p);\ \ \bfY_k(p)=p^{k\beta-1}\bfA_k(p),\ k=1,2,..
\ee

\z where  $\bfA_k$, $k\ge 0$, are analytic functions
 in a neighborhood
of $p=0$ in $\CC$.  
\bigskip

iv) 
If $S_\beta=0$ then $\bfY_k$, $k\ge 0$, are analytic in $\mathcal{W}_1\cup\NN$.

\bigskip

v) The 
analytic continuations of $\bfY_k$ 
along paths in ${\cal R_{\mathrm 1}}$  are in 
$\lloc(\RR^+)$ (their singularities along $\RR^+$
are integrable).
The analytic continuations of the $\bfY_k$
in $\mathcal{R}_1$ can be expressed in terms of each other
through  ``resurgence'' relations
of the type:

\be{resur2}
S_\beta^k\bfY_k=\left(\bfY_0^--\bfY_0^{-^{k-1}+}\right)\circ\tau_{k},\ 
\ \ \ on\ (0,1);\ \ \ (\tau_a:=p\mapsto p-a)
\ee
 
\z relating the higher order series in the trans-series to the first
series and 

\be{resur1}
\bfY_k^{-^m+}=\bfY_k^++\sum_{j=1}^m {k+j\choose k}
S_\beta^j\bfY_{k+j}^+\circ\tau_j
\ee

\ePp
\bigskip

$S_\beta$ is related to the Stokes constant $S$ by

\[ S_\beta=\left\{\begin{array}{l}
\displaystyle\frac{iS}{2\sin(\pi(1-\beta))}\mbox{ for $\beta\ne 1$}\cr
\displaystyle\frac{iS}{2\pi}\mbox{ for $\beta= 1$}\end{array}\right.
 \]

The Borel transformability of the principal series 
$\hatby_0$ has been considered for general systems of differential
equations, allowing for resonances (see \cite{BRBS},\cite{Braaksma}).

Let $\bfY$ be one of the functions  $\bfY_k$ and define, on
$\RR^+\cap\mathcal{R}_1$ the ``balanced average'' of $\bfY$:

\be{defmed}
\bfY^{ba}=\bfY^+ +\sum_{k=1}^{\infty}{2^{-k}}\left(\bfY^{-}
-\bfY^{-^{k-1}+}\right)\heav\circ\tau_k
\ee

\z (${\cal H}$ is Heaviside's
function).
For any value of the argument, only finitely many terms
(\ref{defmed}) are nonzero. Moreover, the balanced average preserves
reality
in the sense that if (\ref{eqor}) is real  and
$\hatby_0$ is real then $\bfY^{ba}$ is real on $\RR^+-\NN$
(and in this case the formula can be
symmetrized
by taking 1/2 of the expression above plus 1/2 of the same expression
with $+$ and $-$ interchanged).
Equation (\ref{defmed}) has the main
features of medianization (cf. \cite{Ecalle}), in particular (unlike
individual analytic continuations, see Appendix~\ref{noncomut}) 
commutes with
convolution (cf. Theorem~\ref{Tecal}). As it will become
clear, the advantage of the
definition (\ref{defmed}) is that $\bfY^{ba}$ is exponentially bounded
at infinity for the functions we are dealing with.

Let again $\hatby$ be one of $\hatby_k$ and $\bfY=\bor\hatby$. We
define:

\ba{deflap}
&&\lap_\phi\bor\hatby:=\lap_\phi\,\bfY=x\mapsto \int_0^{\infty e^{i\phi}}\bfY(p)e^{-px}dp
\ \ \mbox{if $\Phi\ne\RR^+$ }\cr
&&\lap_0\bor\hatby:=\lap_0\,\bfY=x\mapsto \int_0^{\infty}\bfY^{ba}(p)e^{-px}dp\ \ 
\mbox{if $\Phi=\RR^+$}
\end{eqnarray}

The connection between true and formal solutions of the differential equation 
is given in the following theorem:

\Tm{teo2}

i)  
 There is a large enough $b$ 
such that, for $\Re(x)>b$  the Laplace transforms
$\lap_\phi\bfY_k$  exist for all $k\ge 0$ and $\phi\in(-\phi_-,\phi_+)$,
cf. (\ref{defsect}).

For $\phi\in(-\phi_-,\phi_+)$ and  any $C$ the series
\be{gentrsum}
\bfy(x)=(\lap_\phi{\cal B}\hatby_0)(x)+\sum_{k=1}^{\infty}C^ke^{-kx}
(\lap_\phi{\cal B}\hatby_k)(x)\ee

\z is convergent for large enough $x$ in
the right half plane. 

The function $\bfy$ in (\ref{gentrsum})
is a solution of the differential
 equation (\ref{eqor}).

 Furthermore,   for any $k\ge 0$ we have
$\lap_\phi{\cal B}\hatby_k\sim\hatby_k$ in the right half
plane and $\lap_\phi{\cal B}\hatby_k$ is a solution of
 the corresponding equation in 
(\ref{systemform}).

\smallskip
ii) Conversely, given $\phi$,
{\em any} solution of (\ref{eqor}) having $\hatby_0$
as an asymptotic series in the right half plane can be written
in the form (\ref{gentrsum}), for a unique $C$.

\smallskip
iii) The constant $C$, associated in ii) with a given solution $\bfy$ of
(\ref{eqor}),
 depends on the angle $\phi$:

\be{microstokes}
C(\phi)=\left\{\begin{array}{ll}
C(0_+)\ \ \mbox{for }\phi>0\cr
C(0_+)-\frac{1}{2}S_\beta\ \ \mbox{for } \phi=0\cr
C(0_+)-S_\beta\ \ \mbox{for }\phi<0\end{array}\right.
\ee

\z (see also (\ref{caracty0p=1}) ).

\eTm

Note that by (\ref{microstokes})  the change 
in the correspondence (\ref{eqrep1}) occurs when
the Stokes line $\arg x=0$ is crossed. This is a {\em local}
manifestation of the Stokes phenomenon (\cite{Stokes}, \cite{Wasow}, \cite{Sibuya}).

\bigskip

\centerline{**}

\bigskip

Next, we study the correspondence
between the solutions of the differential equations
(\ref{eqor}), (\ref{eqM}), their formal solutions and the
solutions of the inverse Laplace transform 
of these equations, which, in the transformed space,
are convolution equations. 

 With the convolution defined as

\be{defconv}
f*g:=p\mapsto\int_0^p f(s)g(p-s)ds
\ee

\z we have, as is well known, $\lap(f*g)=\lap(f)\lap(g)$,
$\lap(-pf(p))=\lap(f(p))^\prime$. (See Section~\ref{usefulfor}
for a few more useful formulas.) In (\ref{eqor}) we write

\be{Taylor series}
{\bf g}(\xi^{-1},{\bfy})=\sum_{|{\bf l}|\ge 1}{\bf g}_{\bf l}(\xi) {\bf
y}^{\bf l}=\sum_{m\ge 0;|{\bf l}|\ge 1}{\bf g}_{m,\bf l}\xi^m
{\bfy}^{\bf l} \ \ (|\xi|<\xi_0,|\bfy|<y_0)
\ee

\z where by construction ${\bf g}_{0,\bfl}={\bf g}_{1,\bfl}=0$
if $|\bfl|=1$ and the notation
${\bfz}^{\bf l}$  means $z_1^{l_1}\cdot z_n^{l_n}$
and $|{\bf l}|=l_1+..+l_n$. 
The formal inverse Laplace transform of $\bog(x,\bfy(x))$
is given by:

\be{lapdef}
{\cal L}^{-1}\sum_{|\bf l|\ge 1}
{\bf y}(x)^{\bf l}\left(\sum_{m\ge 0}{\bf g}_{m,\bf l}x^{-m}\right)
=\sum_{|\bf l|\ge 1}{{\bf G}}_{\bf l}*{\bf Y}^{*\bf l}+
\sum_{|\bf l|\ge 2}{\bf g}_{0,\bfl}{\bf Y}^{*\bf l}=:{{\calnb}}(\bf Y)
\ee

\z where 

\be{Lapg}
{{\bf G}}_{\bf l}(p)=\sum_{m=1}^{\infty}{\bf g}_{m,\bf l}\frac{p^{m-1}}{m!}\
\ ({{\bf G}}_{1,\bfl}(0)=0\ \mbox{ if }|\bfl|=1)
\ee

\be{defvector}
{{\bf G_{l}}*\bfY^{*\bfl}}\in\CC^n;\ \ ({{\bf G_{l}}*\bfY^{*\bfl}})_j:=
\left({{\bf G}_{\bfl}}\right)_j*Y_1^{*l_1}*..*Y_{n}^{*l_n}
\ee

\z The inverse Laplace transform
of (\ref{eqor}) is the convolution equation:

\be{eqil}
-p{\bf Y}(p)={\bf F}_0(p)-\hat\Lambda {\bf Y}(p)-
\hat B\int_0^p\bfY(s)ds+{\cal  N}({\bf Y})(p)
\ee

\z (see (\ref{lapdef})) where, since $\mathbf{f}_0(x)=O(x^{-2})$,

\be{hypoF}
\bffz(0)=0
\ee

\z By transforming (\ref{systemform}) we get, similarly:

\ba{eqM}
&&(\hat\Lambda-p-k) {\bfY_k}(p)+
\hat B\int_0^p\bfY_k(s)ds-\sum_{j=1}^n\int_0^p(\bfY_k)_j(s)\bfd_j(p-s)ds
=\cr&&\sum_{|\bfl|> 1}\bfdl *
\sum_{\Sigma m=k}*\prod_{i=1}^n*\prod_{j=1}^{l_i}(\bfY_{m_{i,j}})_i=:
\bfR_k(p)\ \ \ (k=1,2,..)
\end{eqnarray}

\z with $\bfdm:=\lapi (\bog^{(\bfm)}(x,\bfy_0)/
\bfm !)$, $\bfd_j:=\lapi(\mathbf{\partial\mathbf{g}(x,\bfy_0)/\partial
{y_{\mathrm j}}})
$ and $*\prod$ standing for the  convolution product.

For a given ray $\Phi$ we consider the equations
(\ref{eqil})  and  (\ref{eqM})
in $\lloc(\Phi)$.
When $\Phi$ is not a Stokes line,
 the description of the solutions is quite simple:

\Pp{convoutside}
i) If $\Phi$ is a ray in $\cald$, then the equation (\ref{eqil})
has a unique solution in $\lloc(\Phi)$, namely $\bfY_0=\bor\hatby_0$.

ii) For any ray in $\mathcal{W}_1$,
the system (\ref{eqil}), (\ref{eqM}) has the general solution
solution $C^k\bfY_k=C^k\bor\hatby_k,\ k\ge 0$.

\ePp

The more interesting case $\Phi=\RR^+$ is dealt with in the following theorem:

\Tm{T0}

i) The general solution in $\lloc(\RR^+)$ of the equation
(\ref{eqil}) can be written in the form:

\be{transconv}
\bfY_C(p)=\sum_{k=0}^{\infty}C^k\bfY_{k}^{ba}(p-k)\heav(p-k)
\ee

\z with $C\in\CC$ arbitrary. 

\smallskip

ii) Near $p=1$, 
 $\bfY_C$ is given by:

\be{struct}
{\bfY_C}(p)=\left\{\begin{array}{cc}
& S_\beta(1-p)^{\beta-1}\bfA(p)+\bfB(p) \ \mbox{for $p<1$}\cr
& C(1-p)^{\beta-1}\bfA(p)+\bfB(p) \ \mbox{for $p>1$}
\end{array}\right. \ (\beta\ne 1)
\ee

$$
{ \bfY_C}(p)=\left\{\begin{array}{cc}
& S_\beta\ln (1-p) \bfA(p)+\bfB(p) \ \mbox{for $p<1$}\cr
& (S_\beta\ln (1-p)+C) \bfA(p)+\bfB(p)\ \mbox{for $p>1$}
\end{array}\right. \ (\beta=1)
$$

\z where $\bfA$ and $ {\bf B}$ extend to analytic functions
in a neighborhood of $p=1$.

\smallskip

iii) With the choice  $\bfY_0=\bfY_{0}^{ba}$, the general solution
of (\ref{eqM}) in $\lloc(\RR^+)$ is $C^k\bfY_{k}^{ba}, k\in\NN$.

\eTm

Comparing (\ref{struct}) with (\ref{caracty0p=1}) we see
that if $S\ne 0$ (which is the generic case) 
 the general solution of
(\ref{eqil}) can be written on the interval $(0,2)$ as a linear combination
of the upper and lower analytic continuations of $\bor\hatby_0$:

\be{combilin}
\bfY_C=\lambda_C\bfY_0^++(1-\lambda_C)\bfY_0^-
\ee

Finally we mention the following result, which shows that
the balanced average, like medianization \cite{Ecalle}, commutes
with convolution.

\Tm{Tecal} If $f$ and $g$ are analytic in ${\cal R}_1$
 then $f*g$ extends analytically in ${\cal R}_1$ and 
furthermore,

\be{resEcalle} 
(f*g)^{ba}=f^{ba}*g^{ba}
\ee
\eTm

As a consequence of the linearity of
the balanced averaging and its commutation  with convolution,  
if   ${\bf \tilde t}_{1,2}$ are the
trans-series
 of the solutions ${\bf f}_{1,2}$
of   differential equations of the type considered
in the present paper  (cf.
(\ref{eqrep1})),
and if
$\lap\bor{\bf \tilde t}_{1,2}={\bf f}_{1,2}$ then

\be{comutsum}
\lap\bor\left({a\bf \tilde t}_1+{b\bf \tilde t}_2\right)=
{a\bf f}_{1}+{b\bf f}_{2}
\ee
Moreover, what is less obvious, we have for the component-wise product 
\be{comutprod}
\lap\bor({\bf \tilde t}_1{\bf \tilde t}_2)={\bf f}_1{\bf f}_2
\ee

\z Borel
summation is in fact an isomorphism between a sub-algebra of trans-series
and a function algebra.

\end{section}

\begin{section}{Proofs and further results}
\begin{subsection}{Outline of the proofs of the main results}

To show the results stated in the previous section, 
we first obtain the general solution in $\lloc$ of the 
convolution system  (\ref{eqM}) in $\cald$ and then, separately, 
on the Stokes line $\RR^+$. We show that 
along a  ray in $\cald$, the solution is unique whereas along
the ray $\RR^+$ there is a one-parameter family of solutions
of the system, branching off at $p=1$.  
We show that any $\lloc$ solution of the system
is (uniformly in $k$) exponentially bounded
at infinity therefore Laplace transformable and  (by the usual properties
of the Laplace transform) these transforms solve (\ref{eqor}).
Conversely, any solution of (\ref{eqor}) with the required
asymptotic properties is inverse Laplace transformable, therefore
it has to be one of the previously obtained solutions of the  equation corresponding
to $k=0$.
We then  study the regularity properties of the solutions of the convolution
equation by local analysis.

Having the complete description
of the family of $\lloc$ solutions we compare
different ways that lead to the same solution and
obtain interesting identities; the identities, together 
with the local properties of the solutions are
instrumental in finding the analytic properties of $\bfY_k$
in $\mathcal{R}_1$.

\emph{Key to the main proofs}. The complete connection with
Equation~(\ref{defmed}) is established in Section~\ref{analyt}.
For the remaining parts: {\emph{Proposition~\ref{basicpp}}}: i)
follows from Proposition~\ref{ansector} and Lemma~\ref{L11};
ii) and  iii) follow from Proposition~\ref{propR1}. The proof
of (\ref{analiticstructn}) is given in Remark~\ref{analytprop}
and iv) is shown in Remark~\ref{analyticY}.
Part v) follows from  
Proposition~\ref{resurgenrel} and Proposition~\ref{propR1}. {\emph{Theorem~\ref{teo2}}}:
i) and ii) follow from Lemma~\ref{lallsol} and  Proposition~\ref{Prest}; iii) is
Equation (\ref{stoketran}). {\emph{Proposition~\ref{convoutside}
}} follows from Proposition~\ref{ansector} and Lemma~\ref{lallsol}.
{\emph{Theorem~\ref{T0}}}: follows from
Proposition~\ref{uniright}, Lemma~\ref{l0,1+eps},
 Proposition~\ref{Uniformnorm}. The proof of
\emph{Theorem~\ref{Tecal}} starts with Proposition~\ref{pmed}
and is continued after it.

\end{subsection}

\begin{subsection}{The convolution equation away
from Stokes rays}

For any star-shaped set ${\cal E}$ in the complex plane containing
the origin (i.e., a region such that the origin can be connected with
any other point in ${\cal E}$ by a straight line segment contained in
${\cal E}$) we denote by $\lray({\cal E})$ the set of functions 
which are locally integrable along each ray in $\cal E$.

\Pp{ansector} 
There is a unique
solution of  (\ref{eqil}) in  $\lray({\cal D})$ (cf. (\ref{defD}))
namely $\bfY_0=\bor\hatby_0$.

This solution is analytic in $\cal D$,  Laplace
transformable along any ray $\Phi$ contained in ${\cal D}$ and
$\lap_\phi \bfY_0$ is a solution of (\ref{eqor}).

\ePp

For the proof we need a few more results.

\Rm{9} There is a constant $K>0$ (independent of 
$p$ and $\bfl$) such that for all $p\in \CC$ and all $\bfl\ge{\bf 0}$ 
\be{exponesti}
|{\bf G}_{\bf l}(p)|_\wedge<K \mu^{|\bf l|} e^{\mu |p|}
\ee
\z for $\mu>\max\{\xi_0^{-1},y_0^{-1}\}$ (cf. (\ref{Taylor series}))
($|{\bf f}|_\wedge:=\max_{1..n}\{|f_1|,..,|f_n|\}$ is an
Euclidean norm; for the definition of $\bf G$ see 
(\ref{Lapg}), (\ref{Taylor series}) and (\ref{eqor})).
\eRm 

{\em Proof.}

\z From the analyticity assumption it follows 
 that

\be{unifesti}
|{\bf g}_{m,\bf l}|_\wedge<\mbox {Const } \mu^{m+|\bfl|}
\ee

\z where the constant is independent on $m$ and $\bfl$.  

\z Then, by (\ref{Lapg}),

$$|{\bf G}_{\bf l}(p)|_\wedge<\mathrm{Const}\ \mu^{|\bfl|+1}\frac{e^{\mu| p|}-1}{\mu |p|}<
\mathrm{Const}\ \mu^{|\bfl|+1}e^{\mu |p|}$$

\Box

Consider the ray segments

\be{defPhi}
\Phi_D=\{\alpha e^{i\phi}:0\le\alpha< D
\}
\ee

\z and along $\Phi_D$ the $\lone$ norm with exponential weight

\be{defl1bd}
\|f\|_{b,\Phi}=\|f\|_{b}:= \int_\Phi e^{-b|p|}|f(p)||dp|=
\int_0^{D}e^{-bt}|f(te^{i\phi})|dt
\ee

\z and the space

$$
\lb(\Phi_D):=\{f:\|f\|_{b}<\infty\}$$

\z (if $D<\infty$, $\lb(\Phi_D)=\lloc(\Phi_D))$. We 
mention the following elementary property:

\Rm{contlap} The Laplace transform $\lap$ is a 
continuous operator from $\lb(\Phi_D)$ to the space of analytic
functions in the half plane $\Re(x)>b$ with the uniform norm.

\eRm
$\Box$

Let   $\cal K\in\CC$ be a bounded domain, 
diam$\,({\cal K})=D<\infty$. On the space of
continuous functions on $\cal K$ we take
the uniform norm with exponential weight:
\be{normunifexp}
\|f\|_u:=D\,\sup_{p\in\cal K}\{|f(p)|e^{-b|p|}\}
\ee

\z (which is equivalent to the usual uniform norm).

Let ${\cal O\subset D}, \ \cal O\ni \rm0$ be a {\em star-shaped,
open set}, diam$({\cal O})=D$ containing a ray segment $\Phi$.
 Let  $\cal A$ be  the space of analytic functions 
$f$ in $\cal O$ such that $f(0)=0$, endowed with the norm (\ref{normunifexp}).

\Pp{1} The spaces $\lb(\Phi_D)$ and $\cal A$ are Banach
algebras with respect to the usual addition of
functions and the convolution (\ref{defconv}). Furthermore

\ba{comutconv}
&&\|f*g\|_{b}\le\|f\|_{b}\|g\|_{b}\ \ (f,g\in\lb(\Phi_D))\cr 
&&\|f*g\|_{u}\le\|f\|_{u}\|g\|_{u}\ \ (f,g\in\mathcal{A})\cr
&&\|f*g\|_{u}\le\|f\|_{u}\|g\|_{b}\ \ (f\in C(\Phi_D),g\in\lb(\Phi_D)\cr&&
\end{eqnarray}

($D=\infty$ is allowed in the first inequality).
\ePp

\z With $F(s):=f(se^{i\phi})$ and $G(s):=g(se^{i\phi})$ we have:

\ba{prfRem}
&&\int_0^{D}dte^{-bt}\left|\int_0^{t}dsF(s)G(t-s)\right|\le
\int_0^{D}dte^{-bt}\int_0^{t}ds|F(s)G(t-s)|=\cr
&&\int_0^{D}\int_{0}^{D-v}e^{-b(u+v)}|F(v)||G(u)|dudv\le\cr&&
\int_0^{D}\int_{0}^{D}e^{-b(u+v)}|F(v)||G(u)|
dudv=\|f\|_{b}\|g\|_{b}\cr
&&
\end{eqnarray}

\z On the other hand,
for $f,g\in\cal A$ we have  $f*g\in\cal A$. Also,

\ba{relativecont}
&&\|f*g\|_{u}=D\sup_{p\in\cal O} e^{-b|p|}\left|\int_0^p
f(s)g(p-s)ds\right|\le\cr &&  D\sup_{p\in\cal O}\int_0^{|p|}
\left|f(te^{i\arg p})e^{-bt}g(p-te^{i\arg p})e^{-b(|p|-t)}\right|dt
\end{eqnarray}
\z which is less than both $\|f\|_{u}\|g\|_{u}$ and 
$\|f\|_{u}\|g\|_{b}$.

\Box

\Rm{R0} For $f$ in $\cal A$ or $f$ in $\lb(\Phi_D)$,

\be{resrem0}
\|f\|_{u,b}\rightarrow 0 \ \ \ {\mbox{as}}\ \ b\rightarrow\infty
\ee

\z where $\|\|_{u,b}$ is either of the $\|\|_u$ or $\|\|_b$ and
$D=\infty$ is allowed in the second case.

\eRm

For $\|\|_b$, Eq. (\ref{resrem0}) is an immediate consequence of the
dominated convergence theorem whereas for $\|\|_u$ it follows from the
definition of $\cal A$.

\Box

\Co{corolla} Let $f$ be continuous along
$\Phi_D$, $D<\infty$  and $g\in
\lb(\Phi_D)$.  Given $\epsilon>0$ there exists a  large enough $b$ and
 $K=K(\epsilon,\Phi_D)$ such that for all $k$

$$ \|f*g^{*k}\|_u<K\,\epsilon^k$$

\eCo

By Remark \ref{R0} we can choose $b=b(\epsilon,{\Phi_D})$ so large that
$\|g\|_{b} <\epsilon$. Then, by Proposition~\ref{1} and
Eq. (\ref{normunifexp}) we have:

$$\left|\int_0^{pe^{i\phi}} f(pe^{i\phi}-s)g^{*k}(s)ds\right|\le
D^{-1}e^{b|p|}\|f\|_u\int_0^{pe^{i\phi}} e^{-b|s|}|g^{*k}(s)||ds|\le$$
\vskip -0.3cm

$$D^{-1}e^{b|p|}\|f\|_u\|g\|_{b}^k
<K\,\epsilon^k$$

\Box

\Rm{2} By (\ref{exponesti}), for any $b>\mu$, and 
$\Phi_D\subset\CC$, $D\le\infty$

\be{est LG}
\|{\bf G}_{\bf l}\|_{b}\le K\mu^{|{\bf l|}}\int_0^{\infty}|dp|e^{|p|(\mu-b)}=K\frac{\mu^{|{\bf l|}}}{b-\mu}
\ee
 \eRm
\z where, to avoid cumbersome notations, we write

\be{defbf}
\bff\in \lb(\Phi_D)\ \mbox{iff }\ 
\||{\bf f}|_\wedge\|_{b}\in \lb(\Phi_D)
\ee

\z (and similarly for other norms of vector functions).

{\em Proof of Proposition~\ref{ansector}.}

 We first show existence and
uniqueness in $\lray({\cal D})$ which amounts to nothing more then
existence and uniqueness along each $\Phi_D\subset\cald$.

Then we show
that for large enough $b$ there exists a unique solution
of (\ref{eqil}) in $\lb(\Phi_{\infty})$. Since
this solution is also in $\lloc(\Phi_{\infty})$ it follows that  
 our (unique) $\lloc$ solution is Laplace
transformable. Analyticity is proven by
finding the solution as a fixed point
of  a contraction with respect to the uniform norm in
 a suitable space of analytic functions.

\Pp{proporem}
i) For $\Phi_D\in\cald$ and large enough $b$, the
operator 

\be{defope}
{\cal N}_1:=\bfY(p)\mapsto (\hat\Lambda-p)^{-1}
\left(\bffz(p)-\hat B\int_0^{p}\bfY(s)ds+\calnb(\bfY)(p)\right)
\ee

\z  is a contraction in a small enough neighborhood
of the origin with respect to $\|\|_{u}$ if $D<\infty$ and
with respect to $\|\|_b$ for $D\le\infty$.

%%%

ii) For $D\le\infty$ the operator $\calnb$ given
formally in (\ref{lapdef})
is continuous in $\lloc(\Phi_D)$. The last sum in (\ref{lapdef})
converges uniformly on compact subsets
of $\Phi_D$. $\calnb(\lloc(\Phi_D))\subset
AC(\Phi_D)$, the absolutely continuous functions on $\Phi_D$. Moreover, if $\mathbf{v}_n\rightarrow\mathbf{v}$
in $\|\|_b$ on $\Phi_D$, $D\le\infty$, then for 
$b'\ge b$ large enough, $\calnb({\bf v}_n)$ exist and converge
in $\|\|_{b'}$ to $\bf v$.

\ePp

The last statements amounts to saying that $\calnb$ is
continuous in the topology of the inductive limit 
of the $\lb$.

{\em Proof}.

\z Since $\hat\Lambda$ and $\hat B$
are constant matrices,

\be{evaln1}
\|{\cal N}_1(\bfY)\|_{u,b}\le
\mbox{Const}(\Phi)\left(\|\bffz\|_{u,b}+
\|\bfY\|_{u,b}\|1\|_{b}
+\|\calnb(\bfY)\|_{u,b}\right)
\ee

\z As both $\|1\|_{b}$ and $\|\bffz\|_{u,b}$ 
are $O(b^{-1})$ for large $b$, the fact that
${\cal N}_1$ maps a small ball into itself
 follows from the following Remark.

\Rm{3} Let $\epsilon>0$ be small enough. Then, there is a $K$ such that for
large $b$ and all  $\bf v$ such that  $\|{\bf v}\|_{u,b}=:\delta<\epsilon$,

\be{est1}
 \|{\calnb}({\bf v})\|_{u,b}\le K\left(\,b^{-1}+\|{\bf v}\|_{u,b}\right)\|{\bf v}\|_{u,b}
 \ee

\eRm
\z By (\ref{unifesti}) and (\ref{est LG}), for large $b$ and some positive constants $C_1,..,C_5$,
\ba{estradball}
&&\|{\calnb}({\bf v})\|_{u,b}\le C_1
\left(\sum_{|{\bf l}|\ge 1}\|{\bf G_l}\|_b\|
{\bf v}\|_{u,b}^{|\bf l|}+\sum_{|{\bf l}|\ge 2}\|
{\bf g_{\rm 0,\bf l}}\|_b\|
{\bf v}\|_{u,b}^{|\bf l|}\right)\cr &&
\le \frac{C_2}{b}\left(\sum_{|{\bf l}|\ge 1}\frac{\mu^{|\bf
l|}}{b-\mu}\delta^{|\bf l|}+\sum_{|{\bf l}|\ge 2}\mu^{|\bf
l|}\delta^{|\bf l|}\right)\le\left( C_2 \sum_{m=1}^\infty+\sum_{m=2}^{\infty}\right)\mu^{m}\delta^{m}\sum_{|{\bf l}|=m} 1\cr&&
\le
\left(\frac{C_4}{b}+\mu\delta\right)\sum_{m=1}^{\infty}\mu^m\delta^m
(m+4)^n\le\left(\frac{C_4}{b}+\mu\delta\right)C_5\delta\cr&&
\end{eqnarray}

\Box

To show that $\calnb_1$ is a contraction we need the
following:

\Rm{4}  
\be{estconvn}
\|\bfh_\bfl\|:=\|({\bf f}+{\bf h})^{*\bfl}-\bff^{*\bfl}\|\le |\bfl|\left(
\|\bff\|+\|\bfh\|\right)^{|\bfl|-1}\|\bfh\|
\ee

\z where $\|\|=\|\|_u$ or $\|\|_b$.
\eRm

This estimate will be useful to us when $\bfh$ is a ``small perturbation''.
The proof of (\ref{estconvn}) 
is a simple induction on $\bfl$, with respect to the 
 lexicographic ordering.
For $|\bfl|=1$, (\ref{estconvn}) is trivial; assume (\ref{estconvn}) holds 
for all $\bfl<\bfl_1$ and that $\bfl_1$ differs from its
predecessor $\bfl_0$ at the  position k (we can take $k=1$), i.e.,
 $(\bfl_1)_1
=1+(\bfl_0)_1$. We have:

\ba{estconvn0}
&&\|({\bf f}+{\bf h})^{*\bfl_1}-\bff^{*\bfl_1}\|=
\|({\bf f}+{\bf h})^{*\bfl_0}*({\bf f}_1+{\bf
h}_1)-\bff^{*\bfl_1}\|=\cr
&&\|(\bff^{*\bfl_0}+\bfh_{\bfl_0})*(f_1+h_1)-\bff^{*\bfl_1}\|=
\|\bff^{*\bfl_0}*h_1+\bfh_{\bfl_0}*f_1+\bfh_{\bfl_0}*h_1\|\le\cr&&\cr
&&
\|\bff\|^{|\bfl_0|}\|\bfh\|+\|\bfh_{\bfl_0}\|\|\bff\|+
\|\bfh_{\bfl_0}\|\|\bfh\|\le\cr&&\cr && \|\bfh\|
\left(\|\bff\|^{|\bfl_0|}+|\bfl_0|(\|\bff\|+\|\bfh\|)^{|\bfl_0|}\right)
\le\cr &&\cr && \|\bfh\|(|\bfl_0|+1)(\|\bff\|+\|\bfh\|)^{|\bfl_0|}
\end{eqnarray}

\Rm{5} For small $\delta$ and large enough $b$, $\calnb_1$ defined in a
ball centered at zero, of radius $\delta$ in the norms $\|\|_{u,b}$ is
contractive.\eRm

By (\ref{evaln1}) and (\ref{est1}) we know
that the ball is mapped into itself for large $b$. Let $\epsilon>0$
be small and let 
$\bf f, h$ be such that
$\|\mathbf{f}\|<\delta-\epsilon, \|\mathbf{h}\|<\epsilon$.
Using (\ref{estconvn}) and the notations 
(\ref{eqil}) (\ref{evaln1}) and $\|\|=\|\|_{u,b}$
we obtain, for some positive
constants $C_1,..,C_4$ and large $b$,

\ba{eqdif}
&&\|\calnb_1(\bff+\bfh)-\calnb_1(\bff)\|\le
C_1\|\left(
\sum_{|\bf l|\ge 2}{\bf g_{\rm 0,\bfl}}\cdot+\sum_{|\bf l|\ge 1}{\bf G}_{\bf l}*\right)\left((\bff+\bfh)^{*\bfl}-\bff^{*\bfl}\right)\|\le\cr&&\cr
&&
C_2\|\bfh\|\left(\sum_{|\bf l|\ge1}\frac{\mu^{|{\bf l|}}}{b-\mu}|
\bfl\|\delta^{|\bfl|-1}+\sum_{|\bf l|\ge 2}|\bf l|\mu^{|{\bf l|}}\delta^{|{\bf l|-1}}\right) <
(C_3b^{-1}+C_4\delta)\|\bfh\|\cr&&
\end{eqnarray}

To finish the proof of Proposition~\ref{proporem} take ${\bf v}\in\cal
 A$.  Given $\epsilon>0$ we can choose $b$ large enough (by Remark
 \ref{R0}) to make $\|{\bf v}\|_u<\epsilon$. Then the sum in the formal
 definition of $\caln$ is convergent in $\cal A$, by (\ref{estradball}).
 Now, if $D<\infty$ $\lloc(\Phi_D)=\lb(\Phi_D)$ for any $b>0$.  If ${\bf
 v}_n\rightarrow{\bf v}$ in $\lb(\Phi_D)$, we choose $\epsilon$ small
 enough, then $b$ large so that $\|{\bf v}\|_b<\epsilon$, and finally
 $n_0$ large so that for $n>n_0$ $\|{\bf v}_n-{\bf v}\|_b<\epsilon$
 (note that $\|\|_{b}$ decreases w.r. to $b$) thus $\|{\bf
 v}_n\|_b<2\epsilon$ and continuity (in $\lb(\Phi_D)$ as well as in
 $\lloc(\Phi_\infty)\equiv \cup_{k\in\Phi_\infty}\lb(0,k))$ follows from
 Remark~\ref{5}. Continuity with respect to the topology of the
 inductive limit of the $\lb$ is proven in the same way.  It is
 straightforward to show that $\calnb(\lloc(\Phi))\subset AC(\Phi)$.

\Box$_{P_{\ref{proporem}}}$

The fact that $\lap_\phi\bfY_0$ is a solution
of (\ref{eqor}) follows from Proposition~\ref{proporem},
from Remark~\ref{contlap} and the  elementary properties of $\lap$
(see also the proof of Proposition~\ref{lemsol}).

Since $\bfY_0(p)$ is analytic for small $p$, $(\lap\bfY_0)(x)$
has an asymptotic series for large $x$, which has to agree
with $\hatby_0$ since $\lap\bfY_0$ solves (\ref{eqor}). This
shows that $\bfY_0=\bor\hatby_0$.

\Box$_{P\ref{ansector}}$

\Rm{0} For any $\delta$ there is a constant $K_2=K_2(\delta,|p|)$
so that for all $\bfl$ we have

%%%
\be{controlYl}
|\bfY_0^{*\bfl}(p)|_\wedge\le K_2\delta^{|\bfl|}
\ee
\eRm

The estimates (\ref{controlYl}) follow immediately 
from analyticity  and
from Corollary~\ref{corolla}. 

\nobreak\Box

\end{subsection}
\end{section}

\begin{subsection}{Behavior of $\bfY_0(p)$ near $p=1$.}

Let $\bfY_0$ be the unique solution in $\lray(\cald)$ of (\ref{eqil}) and let
$\epsilon>0$ be small.  Define

\be{defbfh}
\bfH(p):=\left\{\begin{array}{cc}
\bfY_0(p)\ \ \mbox{for $p\in\cald$\,,$|p|<1-\epsilon$}\cr
0\ \ \ \mbox{otherwise}
\end{array}\right.\ \ \mbox{and}\ \ \bfh(1-p):=\bfY_0(p)-\bfH(p)
\ee

\z In terms of $\bfh$, for real $z=1-p, z<\epsilon$,  the equation  (\ref{eqil}) reads:
\be{eq002}
-(1-z){\bf h}(z)=\mathbf{F}_1(z)-\hat\Lambda \bfh(z)+
\hat B\int_{\epsilon}^z\bfh(s)ds+\calnb({\bf H}+{\bf h})
\ee

\z where 

$${\bf F}_1(1-s):=\bffz(s)-{\hat B}\int_0^{1-\epsilon}\bfH(s)ds
$$

\Pp{P8}

i) For small $\epsilon$, $\bfH^{*\bfl}(1+z)$  extends to an analytic
function in the disk
$\calv:=\{z:|z|<\epsilon\}$. Furthermore,
for any $\delta$ there is an $\epsilon$ and a
 constant $K_1:=K_1(\delta,\epsilon)$ such that for 
$z\in\calv$ the analytic continuation satisfies the estimate

\be{estHl2}
|\bfH^{*\bfl}(1+z)|_\wedge<K_1\delta^{l} 
\ee

\ePp
\smallskip

{\em Proof.}

The case $|\bfl|=1$  is trivial: $\bfH$ itself extends
as the zero analytic function. 
We assume by induction on $|\bfl|$ that Proposition~\ref{P8} is true 
for all $\bfl,\,|\bfl|\le l$ and show that it then holds for 
(e.g.)  $H_1*\bfH^{*\bfl}$, for all $\bfl,\,|\bfl|\le l$.

$\bfH$ is analytic in an $\epsilon$--neighborhood of $[0,1-2\epsilon]$,
 and therefore so is $\bfH^{*\bfl}$.
 Taking first $z\in\RR^+,\ z<\epsilon$, we have

\ba{descomconvo}
&&\int_0^{1-z}H_1(s)\bfH^{*\bfl}(1-z-s)ds=
\int_0^{1-\epsilon}H_1(s)\bfH^{*\bfl}(1-z-s)ds=\cr
&&\int_0^{1/2}H_1(s)\bfH^{*\bfl}(1-z-s)ds+
\int_{1/2}^{1-\epsilon}H_1(s)\bfH^{*\bfl}(1-z-s)ds
\end{eqnarray}

\z The  integral on $[1/2,1-\epsilon]$ is analytic for small $z$,
since the argument of $\bfH^{*\bfl}$ varies 
in an $\epsilon$-neighborhood of $[0,1/2]$; the integral
on $[0,1/2)$ equals

\be{int1}
\int_{1/2-z}^{1-z}H_1(1-z-t)\bfH^{*\bfl}(t)dt=
\left(\int_{1/2-z}^{1/2}+\int_{1/2}^{1-\epsilon}
+\int_{1-\epsilon}^{1-z}\right)H_1(1-z-t)\bfH^{*\bfl}(t)dt
\ee

\z In  (\ref{int1}) the 
integral on $[1/2-z,1/2]$ is  clearly  analytic in
$\calv$, the following one is 
the integral of an analytic function
of the parameter $z$ with respect to the absolutely
continuous measure $\bfH^{*\bfl}dt$
whereas in the last integral,
both  $\bfH^{*\bfl}$ (by induction)  and $H_1$ extend
analytically in $\calv$.

To prove now the the induction step for the estimate (\ref{estHl2}),
 fix $\delta$ small and let:

\be{condchoice}
\eta<\delta;\ M_1:=\max_{|p|<1/2+\epsilon}|\bfH(p)|_\wedge;\ 
M_2(\epsilon):=\max_{0\le x\le 1-\epsilon}|\bfH(p)|_\wedge;\ \ 
\epsilon<\frac{\delta}{4 M_1}
\ee

Let $K_2:=K_2(\eta;\epsilon)$ be large enough so that
(\ref{controlYl}) holds with $\eta$ in place of $\delta$
for {\em real}  $x\in [0,1-\epsilon]$ and {\em also} 
in an $\epsilon$ neighborhood in $\CC$ of the interval
$[0,1/2+2\epsilon]$. We use (\ref{controlYl}) to estimate
the second integral in the decomposition
(\ref{descomconvo}) and the first two 
integrals on the r.h.s. of (\ref{int1}). For the last
integral in (\ref{int1}) we use the induction hypothesis.
If $K_1>2K_2\left(2M_1+M_2\right)$, it follows that $|\bfH^{*\bfl}*H_1|_\wedge$ is bounded by (the terms are
in the order explained above):

\be{boundhl+1}
M_2(\epsilon) K_2\eta^l+M_1 K_2\eta^l+M_1 K_2\eta^l
+(2\epsilon)M_1 K_1 \delta^l< K_1 \delta^{l+1}
\ee

\Box

\Pp{pint}
The  equation (\ref{eq002}) can be written as

\be{formN}
-(1-z){\bf h}(z)={\bf F}(z)-\hat\Lambda \bfh(z)+
\hat B\int_{\epsilon}^z\bfh(s)ds-\sum_{j=1}^n
\int_\epsilon^zh_j(s)\bfd_j(s-z)ds
\ee

\z where 

\be{defD1}
{\bf F}(z):=
\calnb(\bfH)(1-z)+{\bf F}_1(z)
\ee

\be{defderiv}
\bfd_j=
\sum_{|\bfl|\ge 1}l_j\bfgl*\bfH^{*\bar\bfl^j}+
\sum_{|\bfl|\ge 2}l_j{\bf g_{\rm 0,\bf l}}\bfH^{*\bar\bfl^j};\ \bar\bfl^j:=(l_1,l_2,..(l_j-1),..l_n)
\ee

\z (cf. also (\ref{defvector})) extend to analytic functions in $\calv$
(cf. Proposition~\ref{P8}).
Moreover, if $\bfH$ is a vector in $\lb(\RR^+)$ then, for large 
$b$, $\bfd_j\in\lb(\RR^+)$ and the functions
$ {\bf F}(z)$ and ${\bf D}_j$ extend to analytic functions in
$\calv$. 
\ePp

{\em Proof.}

Noting that $(\bfY_0-\bfH)^{*2}(1-z)=0$ for $\epsilon<1/2$ and
 $z\in \calv$
the result is easily obtained by re-expanding
$\calnb(\bfH+\bfh)$ since Proposition \ref{P8}
guarantees the uniform convergence of the series thus obtained.
The proof that ${\bf D}_j\in\lb$ for large $b$
is very similar to the proof of (\ref{eqdif}). The analyticity
properties
follow easily from Proposition~\ref{P8}, since 
the series involved in $\caln(\bfH)$ and ${\bf D}_j$
converge uniformly for $|z|<\epsilon$. %%%%

\Box

Consider again the equation (\ref{formN}).
Let $\ga=\hat\Lambda-(1-z){\hat{1}}$, where
$\hat{1}$ is the identity matrix. By construction
$\ga$ and $\hb$ are  block-diagonal, their first block
is one-dimensional: $\ga_{11}=z$ and $\hb_{11}=\beta$.
We write this as $\ga=z\oplus\gc(z)$ and similarly,
$\hb=\beta\oplus\bc$, where $\gc$ and $\bc$
are $(n-1)\times(n-1)$ matrices. $\gc(z)$ and
$\gc^{-1}(z)$ are analytic in $\calv$.

\Lm{L11}
The function $\bfY_0$ given in Proposition~\ref{ansector}
can be written  in the form

\ba{decomY}
\bfY_0(p)=(1-p)^{\beta-1}\bfa_1(p)+\bfa_2(p)\ \ (\beta\ne 1)
\cr
\bfY_0(p)\phantom{z}=\phantom{z}\ln (1-p)\bfa_1(p)+\bfa_2(p)\ \ (\beta= 1)
\end{eqnarray}

\z for $p$ in the region $(\calv+1)\cap{\cal D}$ ($\calv+1:=\{1+z:z\in\calv\}$)
where $\bfa_1,\ \bfa_2$ are analytic functions
in $\calv+1$ and $(\bfa_1)_j=0$ for
$j>1$.
\eLm

{\em Proof.}

 Let $\bfQ(z):=\int_\epsilon^z\bfh(s)ds$. By Proposition~\ref{ansector},
$\bfQ$ is 
analytic in $\calv\cap(1-{\cal D})$. From (\ref{formN}) we obtain

\be{difeqbfQ}
\zpp\bfQ'(z)-\betapp \bfQ(z)={\bf F}(z)-
\sum_{j=1}^n\int_\epsilon^z\bfd_j(s-z)Q_j'(s)ds
\ee
\z or, after integration by parts in the r.h.s. 
of (\ref{difeqbfQ}), ($\bfd_j(0)=0$, cf. (\ref{defderiv})),

\be{difeqbfQ2}
\zpp\bfQ'(z)-\betapp \bfQ(z)={\bf
F}(z)+\sum_{j=1}^n\int_\epsilon^z\bfd_j'(s-z)Q_j(s)ds
\ee

\z With the notation
$(Q_1,{\bf Q}_\perp):=(Q_1,Q_2,..,Q_n)$
we write the system in the form

\ba{syt1}
&&(z^{-\beta} Q_1(z))'=z^{-\beta-1}\left(F_1(z)+
\sum_{j=1}^n\int_\epsilon^z D_{1j}'(s-z)Q_j(s)ds\right)\cr
&&(e^{\hatc(z)}\Qp)'=e^{\hatc(z)
}\gc(z)^{-1}\left(\FP+\sum_{j=1}^n\int_\epsilon^z
\dpp'(s-z)Q_j(s)ds\right)\cr
&&\hatc(z):=-\int_0^z \gc(s)^{-1}\bc(s)ds
\cr
&&\bfQ(\epsilon)=0
\end{eqnarray}

\z After integration we get:

\ba{eqperp12}
Q_1(z)=R_1(z)+J_1(\bfQ)\cr
\Qp(z)=\Rp(z)+J_\perp(\bfQ)
\end{eqnarray}

\z with
\ba{eqperp2}
&&J_1({\bf Q})=z^{\beta}\int_\epsilon^zt^{-\beta-1}\sum_{j=1}^n\int_\epsilon^tQ_j(s)D_{1j}'(t-s)dsdt\cr
&& 
J_\perp({\bf Q})(z):=e^{-\hatc(z)}
\int_\epsilon^z e^{\hatc(t)}\gc(t)^{-1}
\left(\sum_{j=1}^n\int_\epsilon^z
\dpp'(s-z)Q_j(s)ds\right)dt\cr
&&\Rp(z):=e^{-\hatc(z)}
\int_\epsilon^z e^{\hatc(t)}\gc(t)^{-1}
 {\bf F}_\perp(t)dt
\cr
&&R_1(z)=z^{\beta}\int_\epsilon^zt^{-\beta-1}F_1(t)dt
\phantom{mmmmmmm\frac{F_1(s)-F_1(0)}{s}ds}\ \ \ \ \ \ (\beta\ne 1)\cr
&&R_1(z)=F_1(0)+F'_1(0)z\ln z+
z\int_{\epsilon}^z\frac{F_1(s)-F_1(0)-sF'_1(0)}{s}ds
\ \ (\beta = 1)\cr
&&
\end{eqnarray}

Consider the following space of functions:

\ba{funspace}
{\cal T_\beta}&=&
\Big\{{\bf Q}\ \mbox{analytic in }{\calv\cap(\cald-1)}:
{\bf Q}=z^{\beta}\bfA(z)+\bfB(z)\Big\}\  \mbox{for $\beta\ne 1$ and}\cr
{\cal T}_1&=&\Big\{{\bf Q}\ \mbox{analytic in }{\calv\cap(\cald-1)}:
{\bf Q}=z\ln z\bfA(z)+\bfB(z)\Big\}\
\end{eqnarray}

\z where $\bfA,\bfB$ are analytic in ${\calv}$.
(The decomposition of $\bf Q$ in (\ref{funspace})
is unambiguous since $z^{\beta}$ and
$z\ln z$ are not meromorphic in $\calv$.) 

The norm

\be{normT}
\|{\bf Q}\|=\sup\left\{
|\bfA(z)|_\wedge,|\bfB(z)|_\wedge:z\in\calv\right\}
\ee

\z makes $\cal T_\beta$ a Banach space.

\z For $A(z)$ analytic in $\calv$  the
following elementary identities are useful in what follows:

\ba{convert}
&&\int_\epsilon^z A(s)s^r ds
=Const+z^{r+1}\int _{0}^{1}\!{ A}(zt)t^{r}{dt}=Const+
z^{r+1}Analytic(z)\cr
&&\int_0^z s^r\ln s\,A(s)ds=z^{r+1}\ln z\int _{0}^{1}\!{ A}(zt)t^{r}{dt}+z^{r+1}\int _{0}^{1
}\!{ A}(zt)t^{r}\ln t{dt}\cr
&&
\end{eqnarray}

\z where the second equality is obtained  by differentiating
 with respect to $r$  the first equality.

Using
(\ref{convert}) it is straightforward to check that the r.h.s.
of (\ref{eqperp12}) extends to a linear inhomogeneous
operator
on $\cal T_\beta$ with image in $\cal T_\beta$ and that
the norm of $J$ is $O(\epsilon)$ for small $\epsilon$.
For instance, one of the terms in $J$ for $\beta=1$,

\ba{arr1}
&&z\int_0^z t^{-2}
\int_0^t s\ln s\,A(s)D'(t-s)ds=\cr &&z^{2}\ln z\int _{0}^{1}
\int _{0}^{1}\sigma A(z\tau\sigma)D'
(z\tau-z\tau
\sigma){d\sigma}{d\tau}+
\cr &&
z^{2}\int _{0}^{1}d\tau
\int _{0}^{1}d\sigma\sigma(\ln\tau+\ln\sigma)
A(z\tau\sigma)D'(z\tau-z\tau\sigma)
\cr &&
\end{eqnarray}

\z  manifestly in $\cal T_\beta$ if $A$ is
analytic in $\calv$. Comparing with (\ref{funspace}), the extra power 
of $z$
accounts for a norm $O(\epsilon)$ for this term.

Therefore, in (\ref{syt1}) $(1-J)$ is invertible and
the  solution $\bf Q\in\cal T_\beta\subset\cal
L\left(D\right)$. In view of  the the uniqueness of $\bfY_0$
(cf. Proposition
\ref{ansector}), the rest of the proof of Lemma
{\ref{L11}} is immediate.

\Box

\end{subsection}

\begin{subsection}{The solutions of (\ref{eqil}) on $[0,1+\epsilon]$}

Let $\bfY_0$ be the solution given by Proposition
\ref{ansector}, take $\epsilon$
small enough 
and denote by $\cal O_\epsilon$ a neighborhood in $\CC$ of width
$\epsilon$  of the interval $[0,1+\epsilon]$.

\Rm{uniformlb}.
$\bfY_0\in \lone({\cal O}_\epsilon)$. As $\phi\rightarrow\pm 0$,
$\bfY_0(pe^{i\phi})\rightarrow\bfY_0^\pm(p)$ in the
sense of $\lone([0,1+\epsilon])$ and also in the sense of
 pointwise convergence for 
$p\ne 1$, where

\ba{defsolpm}
&&\bfY_0^\pm:=\left\{\begin{array}{cc}
\bfY_0(p)
\phantom{\pm0i)^{\beta-1}\bfa_1(p)+\bfa_2(p)}&{p<1}\cr
(1-p\pm0i)^{\beta-1}\bfa_1(p)+\bfa_2(p)\ \  &{p>1}
\cr
\end{array}\ \ (\beta\ne 1)
\right.\cr
&&\cr
&&\bfY_0^\pm:=\left\{\begin{array}{cc}
\bfY_0(p)
\phantom{1-p\pm 0i)\bfa_1(p)+\bfa_2(p)}&{p<1}\cr
\ln(1-p\pm 0i)\bfa_1(p)+\bfa_2(p)&{p>1}
\cr
\end{array}\ \ (\beta= 1)
\right.
\end{eqnarray}

\z Moreover, $\bfY_0^{\pm}$ are $\lloc$ 
solutions of the convolution equation
(\ref{eqil}) on the interval $[0,1+\epsilon]$.
\eRm

The proof is immediate from Lemma~\ref{L11} and Proposition~\ref{proporem}.

\Box

\Pp{lincomb}
For any $\lambda\in\CC$ the combination
$\bfY_\lambda=\lambda\bfY_0^++(1-\lambda)\bfY_0^-$
is a solution of (\ref{eqil}) on $[0,1+\epsilon]$.

\ePp

{\em Proof. } For $p\in [0,1)\cup(1,1+\epsilon]
$ let ${\bf y}_\lambda(p):=\bfY_\lambda-\bfH(p)$. 
Since ${\bf y}_\lambda^{*2}=0$ the equation
(\ref{eqil}) is actually linear in ${\bf y}_\lambda$
(compare with (\ref{formN})). 

\Box

\centerline{*}

Note: We consider the application $\mathcal{Y}:=\bfy_0\mapsto
\bfY_\lambda$ and require that it is compatible with complex conjugation
of functions $\mathcal{Y}({\bfy_0}^*)= ({\mathcal{Y}({\bfy_0})})^*$
where $F^*(z):=\overline{F(\overline{z})}$. We get $\Re\,\lambda=1/2$.
It is natural to choose $\lambda=1/2$ to make the linear combination a
true average.  This choice corresponds, on $[0,1+\epsilon]$, to the
balanced averaging (\ref{defmed}).

\centerline{*}

\Rm{Rright} For any $\delta>0$  there
is a constant $C(\delta)$ such that for large $b$

\be{estimcv}
\|(\bfY_0^{ba})^{*\bfl}\|_u<C(\delta)\delta^{|\bfl|}
\ \ \forall\,\bfl\ \mbox{with}\ |\bfl|>1
\ee

\z ($\|\|_u$ is taken on the interval $[0,1+\epsilon]$).

\eRm

Without loss of generality, assume that $l_1>1$. Using the notation
(\ref{defderiv}),

\ba{equsplit}
&&\left\|\int_0^p(\bfY_0)_1^{ba}(s)
(\bfY_0^{ba})^{*\bar\bfl^1}(p-s)ds\right\|_u\le\cr
&&\left\|\int_0^{\frac{p}{2}}(\bfY_0^{ba})_1(s)(\bfY_0^{ba})^{*\bar\bfl^1}(p-s)ds\right\|_{u_2}+\left\|
\int_0^{\frac{p}{2}}(\bfY_0)_1(p-s)(\bfY_0^{ba})^{*\bar\bfl^1}(s)ds\right\|_{u_2}
\cr&&
\end{eqnarray}

\z ($\|\|_{u_2}$ refers to the interval $p\in[0,1/2+\epsilon/2]$.)
The first $u_2$ norm can be estimated directly using
Corollary~\ref{corolla}
whereas we majorize the second one by

$$\|(\bfY^{ba}_0)_1\|_b\|(\bfY_0^{ba})^{*\bar\bfl^1}(x)\|_
{u_2}$$

\z and apply Corollary~\ref{corolla} to it for $|\bfl|>2$
(if $|\bfl|=2$ simply observe that
$(\bfY_0^{ba})^{*\bfl}$ is analytic on $[0,1/2+\epsilon/2]$).

\Box

\Lm{l0,1+eps} 

The set of all
solutions of (\ref{eqil})
in $\lloc([0,1+\epsilon])$ is parameterized by a complex
constant $C$ and is given by 

\be{solgen}
\bfY_0(p)=\left\{\begin{array}{cc}
\bfY_0^{ba}(p)\phantom{+Cz^{\beta-1}\bfA(p)+1}\ \mbox{for $p\in[0,1)$}\cr
\bfY_0^{ba}(p)+C(p-1)^{\beta-1}\bfA(p) \ \mbox{for $ p\in(1,1+\epsilon]$}
\end{array}\right.
\ee

\z for $\beta\ne 1$ or, for $\beta=1$,

$$
\bfY_0(p)=\left\{\begin{array}{cc}
\bfY_0^{ba}(p)\phantom{C(p-1)\bfA(p)}\ \mbox{for $p\in[0,1)$}\cr
\bfY_0^{ba}(p)+C(p-1)\bfA(p)\ \mbox{for $ p\in(1,1+\epsilon]$}
\end{array}\right.
\eqno(\ref{solgen})'$$

\z  where $\bfA$ extend 
analytically in a neighborhood of $p=1$. 

Different values of $C$ correspond to
different solutions. 

This result remains true if $\bfY_0^{ba}$  is replaced
by any other combination $\bfY_\lambda:=\lambda\bfY_0^++(1-\lambda)\bfY_0^-$, $\lambda\in\CC$.

\eLm

{\em Proof.}

We look for solutions of (\ref{eqil}) in the form 
\be{defYY}
\bfY^{ba}(p)+\bfh(p-1)\ee

\z By Lemma~\ref{L11} ,
$\bfh(p-1)=0$ for $p<1$. Note that

\be{decomp}
\calnb(\bfY_0^{ba}\circ\tau_{-1}+\bfh)(z)=\calnb(\bfY_0^{ba})(1+z)+
\sum_{j=1}^n\int_0^zh_j(s)\bfd_j(z-s)ds
\ee

\z where the $\bfd_j$ are given in (\ref{defderiv}), and
by Remark~\ref{estimcv} all the infinite sums involved
are uniformly convergent. For $z<\epsilon$ (\ref{eqil})
translates to  
(compare with (\ref{formN})):

\be{eqilh}
-(1+z){\bf h}(z)=-\hat\Lambda {\bfh}(z)-
\hat B\int_0^z\bfh(s)ds+\sum_{j=1}^n\int_0^zh_j(s)\bfd_j(z-s)ds
\ee

\z Let 

\be{defQQ}
\bfQ(z):=\int_0^z\bfh(s)ds
\ee

\z As we
are looking for solutions
$\bfh\in\lone$, we have $\bfQ\in AC[0,\epsilon]$ and
$\bfQ(0)=0$.

\z Following the same steps as in the proof of Lemma
\ref{L11} we get the system of equations:

\ba{sytright}
&&(z^{-\beta} Q_1(z))'=z^{-\beta-1}
\sum_{j=1}^n\int_0^z D_{1j}'(z-s)Q_j(s)ds\cr
&&(e^{\hatc(z)}\Qp)'=e^{\hatc(z)
}\gc(z)^{-1}\sum_{j=1}^n\int_0^z
\dpp'(z-s)Q_j(s)ds\cr
&&\hatc(z):=-\int_0^z \gc(s)^{-1}\bc(s)ds
\cr
&&\bfQ(0)=0
\end{eqnarray}

\z which by integration gives

\ba{eqperpright12}
(\hat 1+J)\bfQ(z)=C{\bf R}(z)
\end{eqnarray}

\z where  $C\in\CC$ and
\ba{eqperpright}
&&(J({\bf Q}))_1(z)=z^{\beta}\int_0^zt^{-\beta-1}\sum_{j=1}^n\int_0^tQ_j(s)D_{1j}'(t-s)dsdt\cr
&& 
J({\bf Q})_\perp(z):=e^{-\hatc(z)}
\int_0^z e^{\hatc(t)}\gc(t)^{-1}
\left(\sum_{j=1}^n\int_0^z
\dpp'(z-s)Q_j(s)ds\right)dt\cr
&&{\bf R}_{\perp}=0\cr
&&R_1(z)=z^{\beta}
\cr
&&
\end{eqnarray}

First we note the presence of an arbitrary constant $C$
in (\ref{eqperpright12}) (Unlike in
Lemma \ref{L11} when the initial condition, given at $z=\epsilon$ was
determining the integration constant, now the initial condition
$\bfQ(0)=0$ is satisfied for all $C$).

For small $\epsilon$ 
the norm of the operator $J$ defined 
on $AC[0,\epsilon]$ is $O(\epsilon)$, as in the proof
of Lemma~\ref{L11}. Given
$C$ the solution of the system (\ref{sytright})
is unique and can be written as 

\be{eqneuman}
\bfQ=C \bfQ_0;\ \ \bfQ_0:=(\hat 1+J)^{-1}{\bf R}\ne 0
\ee

It remains to find the analytic structure of ${\bf Q}_0$. We now
introduce the space

\ba{funspace2}
{\cal T}=
\left\{{\bf Q}:[0,\epsilon)\mapsto\CC^n:
{\bf Q}=z^{\beta}\bfA(z)\right\}
\end{eqnarray}

\z where $\bfA(z)$ extends to an analytic
function
in $\calv$. With the norm (\ref{normT}) (with ${\bf B\equiv 0}$),
$\cal T$ is a Banach space.
As in the proof of Lemma \ref{L11}
the operator $J$ extends naturally to $\cal T$
where it has a norm $O(\epsilon)$ for small $\epsilon$.
It follows immediately that

\be{rezQ0}
{\bf Q}_0\in\cal T
\ee

The formulas
(\ref{solgen}), (\ref{solgen}') follow from
(\ref{defYY}) and  (\ref{defQQ}).

\Box

\Rm{Rnew}
If $S_\beta\ne 0$ (cf. Lemma \ref{L11})
then the \emph{general} solution of (\ref{eqil})
is given by

\be{lincombe}
\bfY_0(p)=(1-\lambda)\bfY_0^+(p)+\lambda\bfY_0^-(p)
\ee

\z with $\lambda\in\CC$.

\eRm

Indeed, if $\bfa_1\not\equiv 0$ (cf. Lemma \ref {L11}) we get at least
two distinct solutions of (\ref{eqperpright12}) (i.e., two distinct
values of $C$) by taking different values of $\lambda$ in
(\ref{lincombe}).  The remark  follows from (\ref{rezQ0})
(\ref{funspace2}) and Lemma~\ref{l0,1+eps}..

\Box

\end{subsection}

\begin{subsection}{The solutions of (\ref{eqil}) on $[0,\infty)$}

In this section we show that the leading asymptotic
behavior of $\bfY_p$ 
as $p\rightarrow 1_+$ determines
a unique solution of (\ref{eqil}) in $\lloc(\RR^+)$. Furthermore,
any $\lloc$ solution of (\ref{eqil})
is exponentially bounded at infinity and thus Laplace transformable. 
We also study some properties of these solutions 
and of their Laplace transforms.

 Let  ${\bf H}$ be a solution of (\ref{eqil})
on an interval $[0,1+\epsilon]$, which we extend to $\RR^+$
letting $\bfH(p)=0$ for $p>1+\epsilon$. For a large enough $b$,
define

\be{spacesol}
{\cal S}_{\bf H}:=\{f\in\lloc([0,\infty)):f(p)={\bf H}(p) \mbox{ on } [0,1+\epsilon]\}
\ee

\z and

\be{defcals0}
{\cal S}_0:=\{f\in\lloc([0,\infty)):f(p)=0 \mbox{ on } [0,1+\epsilon]\}
\ee

\z We extend ${\bf H}$ to $\RR^+$ by
putting ${\bf H}(p)=0$ for $p>1+\epsilon$; for
$p\ge 1+\epsilon$ (\ref{eqil}) reads:

\be{eq2}
-p({\bf H}+{\bf h})=F_0-\hat\Lambda ({\bf H}+{\bf h})-\hat{B}\int_0^p({\bf H}+{\bf h})(s)ds
+{\cal  N}({\bf H}+{\bf h})
\ee

\z with $\bfh\in {\cal S}_0$, or

\be{eqi3}
{\bf h}=-{\bf H}+(\hat\Lambda-p)^{-1}\left(F_0-\hat{B}\int_0^p({\bf H}+{\bf h})(s)ds
+{\caln}({\bf H}+{\bf h})\right):={\cal M}({\bf h})
\ee

\z For small $\phi_0>0$ and $0\le\rho_1<\rho_2\le\infty$,
consider the truncated sectors

\be{truncsect}
S^{\pm}_{(\rho_1,\rho_2)}:=\{z:z=\rho e^{\pm i\phi},\rho_1<\rho<\rho_2;\
0\le\phi<\phi_0\}
\ee

\z and the spaces of functions analytic in $S^{\pm}_{(\rho_1,\rho_2)}$
and continuous in its closure:
\be{defT}
{\cal T}^{\pm}_{\rho_1,\rho_2}
=\left\{f:f\in C(\overline{S_{(\rho_1,\rho_2)}});\ f\ \mbox{analytic
in}\ 
S^\pm_{(\rho_1,\rho_2)}\right\}
\ee

\z  which are Banach spaces with respect to $\|\|_u$ on compact
subsets of $\overline{S_{(\rho_1,\rho_2)}}$.

\Pp{uniright}
i) Given $\bfH$, the equation (\ref{eqi3}) has a unique solution in
$\lloc[1+\epsilon,\infty)$. For large 
$b$, this solution is in
$\lb([1+\epsilon,\infty))$  and thus Laplace
transformable.

ii) Let $\bfY_0$ be the solution  defined in Proposition \ref{ansector}.
Then

\be{regularity}
\bfY_0^{\pm}(p):=\lim_{\phi\rightarrow \pm 0} \bfY_0(pe^{i\phi})\in
C(\RR^+\backslash\{1\})\cap\lloc(\RR^+)
\ee

\z (and the limit exists pointwise on $\RR^+\backslash\{1\}$
and in $\lloc(\RR^+)$.) 

Furthermore, $\bfY_0^{\pm}$ are particular solutions of (\ref{eqil})
and

\ba{descont}
\bfY_0^\pm(p)&=&(1-p)^{\beta-1}\bfa^\pm(p)+\bfa_1^\pm(p)\ \ (\beta\ne 1)
\cr
\bfY_0^\pm(p)&=&\ln (1-p)\,\bfa^\pm(p)+\bfa_1^\pm(p)\ \ (\beta= 1)
\end{eqnarray}

\z where $\bfa^\pm$ and $\bfa_1^\pm$ are in ${\cal T}^\pm_{0,\infty}$.

\ePp

{\em Proof}

Note first that by Proposition \ref{proporem},
${\cal M}$ (eq. (\ref{eqi3}))
is well defined on ${\cal S}_0$, (eq.(\ref{defcals0})).
Moreover, since $\bfH$ is a solution of (\ref{eqil})
on $[0,1+\epsilon)$, we have, for $\bfh_0\in{\cal S}_0$, 
${\cal M}(\bfh)=0$ a.e.
on $[0,1+\epsilon)$, i.e.,

$${\cal M}({\cal S}_0)\subset{\cal S}_0$$

\Rm{new2}
For large $b$, ${\cal M}$ is a contraction 
in a small neighborhood of the origin in $\|\|_{u,b}$.

\eRm
Indeed, $\sup\{\|(\hat\Lambda-p)^{-1}\|_{\CC^n\mapsto\CC^n}:p\ge 1+\epsilon\}
=O(\epsilon^{-1})$ so that

\be{relativenorm}
\|\calm(\bfh_1)-\calm(\bfh_2)\|_{u,b}\le\frac{\mbox{Const}}{\epsilon}\|\n(\bff+\bfh)-\n(\bff)\|_{u,b}
\ee

\z The rest follows from (\ref{eqdif})
---Proposition~\ref{proporem} and Remark~\ref{R0} applied to $\bfH$.

\Box

The existence of a solution of (\ref{eqi3}) in ${\cal
 S}_0\cap\lb([0,\infty))$ for large enough $b$ is now immediate.

Uniqueness in $\lloc$ is tantamount to uniqueness in
 $\lone([1+\epsilon,K])=\lb([1+\epsilon,K]$, for all
 $K-1-\epsilon\in\RR^+$.  Now, assuming $\calm$ had two fixed points in
 $\lb([1+\epsilon,K])$, by Remark \ref{R0}, we can choose $b$ large
 enough so that these solutions have arbitrarily small norm, in
 contradiction with Remark~\ref{new2}.

$ii)$. For $p<1, \bfY_0^{\pm}(p)=\bfY_0(p)$.
For $p\in(1,1+\epsilon)$ the result follows 
from Lemma \ref{L11}. Noting  that (in view of the estimate (\ref{estradball}))
$\calm({\cal T^\pm}_{1+\epsilon,\infty})\subset {\cal
T^\pm}_{1+\epsilon,\infty}$, the
rest of the proof follows from the Remark~\ref{new2}
and  Lemma \ref{L11}.

\Box

\Pp{lemsol} There is a one parameter
family of solutions of equation (\ref{eqil}) in $\lloc[0,\infty)$,
branching off at $p=1$ and in a neighborhood of $p=1$ all solutions
are of the form (\ref{solgen}), (\ref{solgen}').
The general solution of (\ref{eqil})
is Laplace transformable for large $b$ 
and the Laplace transform is a solution
of the original differential equation in the half-space
$\Re( x)>b$.

\ePp

Note: 
As of now, the correspondence (\ref{solgen}), (\ref{solgen}')
with the balanced average (\ref{defmed})
is proven only near $p=1$; the complete correspondence is
established  in Section~\ref{analyt}.

{\em Proof.}

Let ${\bf Y}$ be any solution of (\ref{eqil}).
By  Lemma \ref{l0,1+eps} and 
Proposition~\ref{uniright},  $b$ large
 implies that $\bfY\in\lb([0,\infty))$ (thus $\lap\bfY$ exists), 
that $\|\bfY\|_b$ is small and, in particular, that the sum 
defining $\calnb$ in (\ref{lapdef}) is convergent in
$\lb(\RR^+)$.

By Remark~\ref{contlap},

\ba{comutlaplac}
&&{\cal L}\left(\sum_{|\bf l|\ge1}
{\bf G}_{\bf l}*\bfY^{*\bfl}+\sum_{|\bfl|\ge 2}{\bf g_{\rm 0,\bf l}}
\bfY^{*\bfl}\right)(x)=\cr&&
\sum_{|\bf l|\ge1}
({\cal L}{\bf G}_{\bf l})({\cal L}\bfY)^{\bfl}(x)+
\sum_{|\bfl|\ge 2}{\bf g_{\rm 0,\bf l}}\left(\lap\bfY\right)^{\bfl}=
\sum_{|\bf l|\ge1}
{\bf g}_{\bf l}(x){\bf y}^{\bfl}(x)={\bf g}(x,{\bf y}(x))\cr&&
\end{eqnarray}

\z (and ${\bf g}(x,\bfy(x))$ is analytic for $\Re (x)> b$).  The rest
is straightforward.

\Box

\end{subsection}

\begin{subsection}{Correspondence with formal solutions}
\label{cor-for}

\z Finally we consider  formal solutions 
{\em for large argument} 
of the differential equation, in the differential
algebra generated by formal power series
(in decreasing powers of the large variable)
 and (decreasing) exponentials,
i.e. solutions
as formal asymptotic expansions. The theory
of formal solutions is classical (\cite{Fabry}, \cite{Cope} \cite{Ritt});
see also \cite{Ecalle} for a vast and
very interesting generalization.
We only sketch the facts that are 
relevant to us.

The simplest 
formal solution of (\ref{eqor}) is an asymptotic series $\hatby_0$.

$$\hatby_0=\sum_{m=2}^\infty\frac{\bfy_{0,m}}{x^m}$$

\z In view of the invertibility of $\hat \Lambda$, the coefficients
$\{{\bfy}_{0,m}\}_{m\in\NN}\subset\CC^n$ can be determined uniquely by
expanding in (\ref{eqor}) in powers of $1/x$ and equating  the
coefficients of the $x^{-m},m\ge 2$.  The series $\hatby_0$ is
generically divergent.
 
Since we expect an $n-parameter$ family
of solutions, we look for further 
solutions as perturbations of $\hatby_0$.
Because of the uniqueness
of $\hatby_0$ a perturbation must be  smaller than all powers
of $x^{-1}$ i.e., ``beyond all orders'' of $\hatby_0$.

Taking $\hatby=\hatby_0+\hatby_1$ we get, 
to the lowest order of approximation,
$\hatby_1'=-\hat\Lambda\hatby_1$. The solutions to this approximate 
equation
are linear combinations of $e^{-\lambda x}$,
$\lambda\in$spec$\hat\Lambda$. 
 We only consider  solutions $\hatby_1$ that are (formally)
small perturbations of $\hatby_0$ in the
half-plane $\Re(x)>0$; this condition selects out the eigenvalue
$\lambda=1$.

Continuing the perturbative procedure until we reach
a formal solution of (\ref{eqor}), we end up with
an exponential series

\be{forms1}
\hatby=\hatby_0+\sum_{k=1}^\infty e^{-kx}\hatby_k  
\ee

\z where $\hatby_k$ are formal power series.
Substituting (\ref{forms1}) into (\ref{eqor})
and using the fact that $\hatby_0$ is already
a formal solution we get for $\hatby_k,k\ge 1$:

\ba{partsum}
&&\sum_{k=1}^\infty e^{-kx}\left[\hatby'_k-
\left(k-\hat\Lambda-\frac{1}{x}\hat B+
{\partial}{\bf g}(x,\hatby_0)\right)
\hatby_k\right]=\cr
&&
\sum_{|\bfl|> 1}
\frac{\bogl(x,\hatby_0)}{\bfl!}
\left(\sum_{k=1}^\infty e^{-kx}\hatby_k\right)^{\bfl}=\cr
&&\sum_{k=2}^\infty e^{-kx}\sum_{|\bfl|> 1}\frac{\bogl(x,\hatby_0)}{\bfl!}
\sum_{\scriptscriptstyle\Sigma m=
}\prod_{i=1}^n\prod_{j=1}^{l_i}
(\hatby_{m_{i,j}})_i\cr
&&
\end{eqnarray}

\z Equating the coefficients of $e^{-kx},
k\ge 0$ we get the system (\ref{systemform}).

By assumption, $\hat\Lambda-1$ has a one-dimensional null-space. Thus, by
(\ref{systemform}),  $\hatby_1$ has the freedom of an arbitrary
multiplicative constant. We make a definite choice of $\hatby_1$ by
requiring that the first component of the coefficient of the leading
power of $x$ is one.

Still by assumption, for $k\ne 1$ 
$\hat\Lambda-k$ is
invertible, so that, taking $C=1$,
all $\hatby_{k}$, $k\ge 1$,  are uniquely determined.
Letting $C$ be arbitrary
we get instead  $C\hatby_1$ for $k=1$, 
 $C^2\hatby_2$ for $k=2$ (because of the condition $\sum m=2$),
etc, so that the general formal
solution of type (\ref{forms1}) is 

$$\hatby=\hatby_0+\sum_{k=1}^\infty C^ke^{-kx}\hatby_{k}$$

The existence of formal exponential solutions
has been considered in \cite{Fabry},
\cite{Cope}, \cite{Iwano} and a very comprehensive
theory can be found in Ecalle   \cite{Ecalle-book}, 
\cite{Ecalle}, \cite{Ecalle2}.

The following proposition is a classical result
and is a specialization
of general theorems (see \cite{Iwano}).

\Pp{Pclassical}
 There is exactly a one parameter family of solutions of (\ref{eqor})
having the asymptotic behavior described by 
$\hatby_0$ in the half-plane $\Re(x)>0$.

\ePp

{\em Proof.}  Any solution
with the properties stated in Proposition~\ref{Pclassical} is inverse Laplace
transformable and its inverse Laplace transform has to be one of the
$\lloc$ solutions of the convolution equation (\ref{eqil}).
The rest of the proof follows from Proposition~\ref{lemsol}.

\Box

\Pp{estimd} Let $\bfY$ be any $\lloc(\RR^+)$ solution of (\ref{eqil}).
 For large $b$ and some $\nu>0$ the coefficients
$\bfdm$  in (\ref{eqM}) are bounded by
$$|\bfdm(p)|_\wedge\le e^{\mu p}\nu^{|\bfm|}$$
\ePp

Note that $\lapi (\bog^{(\bfm)}(x,\bfy)/
\bfm !)$
is the coefficient of $\bfZ^{*\bfm}$ in the expansion of
$\calnb(\bfY+\bfZ)$ in convolution powers of $Z$  (\ref{lapdef}):

\ba{exprgl}&&\left(\left(\sum_{|\bf l|\ge 2}{\bf g_{\rm 0,\bf l}}\cdot+\sum_{|\bf l|\ge 1}{{\bf G}}_{\bf l}
*\right)({\bf Y}+\bfZ)^{*\bf
l}\right)_{\bfZ^{*\bfm}}=\cr&&
\left(\left(\sum_{|\bf l|\ge 2}{\bf g_{\rm 0,\bf l}}\cdot+\sum_{|\bf l|\ge 1}{{\bf G}}_{\bf l}
*\right)
\sum_{0\le\bfk\le\bfl}{\bfl\choose\bfk}
\bfZ^{*\bfk}{\bf Y}^{*(\bfl-\bfk)}\right)_{\bfZ^{*\bfm}}=\cr&&
\left(\sum_{|\bf l|\ge 2}{\bf g_{\rm 0,\bf l}}\cdot+\sum_{|\bf l|\ge 1}{{\bf G}}_{\bf l}
*
\right)\sum_{\bf l\ge\bf m}{
{\bfl\choose\bfm}{\bf G}_{\bf l}}*
{\bf Y}^{*(\bfl-\bfm)}
\end{eqnarray}

\z ($\bfm$ is fixed) where $\bf l\ge m$ means $l_i\ge m_i, i=1..n$
and
${{\bfl}\choose{\bfk}}:=\prod_{i=1}^n{{l_i}\choose{k_i}}$.

Let $\epsilon$ be small and $b$ large so that
$\|\bfY\|_b<\epsilon$.
Then, for some constant $K$, estimate (cf. (\ref{exponesti}))

$$
\left|\left(\sum_{II}{\bf g_{\rm 0,\bf l}}\cdot+\sum_{I}{{\bf G}}_{\bf l}
*
\right){
{\bfl\choose\bfm}{\bf G}_{\bf l}}*
{\bf Y}^{*(\bfl-\bfm)}\right|_\wedge
\le
\sum_{I}
{Ke^{\mu|p|}(\mu\epsilon)^{|\bfl-\bfm|}}{\bfl\choose\bfm}=$$

\be{estnorm}
{\epsilon^{-|\bfm|}Ke^{\mu|p|}}
\prod_{i=1}^n\sum_{l_i\ge m_i}
{l_i\choose m_i}(\mu\epsilon)^{l_i}={K}
\frac{e^{\mu|p|}\mu^{|\bfm|}}{(1-\epsilon\mu)^{|\bfm|+n}}<{e^{\mu|p|}\nu^{|\bfm|}}
\ee

\z (where $I(II)\equiv\{|\bfl|\ge 1(2);\bfl\ge\bfm\})$ for  large enough $\nu$.

\Box

For $k=1$, $\bfR_1=0$ and equation (\ref{eqM}) is (\ref{eqilh}) (with $p
\leftrightarrow z$)
but now on the whole line $\RR^+$. For
small $z$ the solution is given by (\ref{eqneuman})
(note that $\bfd_1={\bf d}_{(1,0,..,0)}$
and so on)
and depends on the free constant $C$ (\ref{eqneuman}). We choose 
a value for $C$  (the values of $\bfY_1$
 on $[0,\epsilon]$ are then determined)
and we write  the equation of $\bfY_1$ for $p\ge\epsilon$:

\ba{equextend}
&&(\hat\Lambda-1-p) {\bfY_1}(p)+
\hat B\int_\epsilon^p\bfY_1(s)ds-\sum_{j=1}^n\int_\epsilon^p(\bfY_1)_j(s)\bfd_j(p-s)ds
=
\cr&&\bfR(p):=\int_0^\epsilon\bfY_1(s)ds+\sum_{j=1}^n
\int_0^\epsilon (\bfY_1)_j(s)\bfd_j(p-s)ds
\end{eqnarray}

\z ($\bfR$ only depends on the values of $\bfY_1(p)$ on $[0,\epsilon]$).
We write 

\be{eqabstract}(1+J_1)\bfY_1=\hat Q_1^{-1}\bfR\ee

\z with $Q_1=1-\hat{\Lambda}+p$. The operator $J_1$ is
defined by $(J_1\bfY_1)(p):=0$ for $p<\epsilon$, while, for $p>\epsilon$,
$$(J_1\bfY_1)(p):=Q_1^{-1}\left(\hat
B\int_\epsilon^p\bfY_1(s)ds-\sum_{j=1}^n\int_\epsilon^p(\bfY_1)_j(s)\bfd_j(p-s)ds\right)$$

\z By Proposition~\ref{pint}, (\ref{comutconv}) and Remark~\ref{R0},
 noting that $\sup_{p>\epsilon}\|Q_1^{-1}\|=O(\epsilon^{-1})$, b we find
 that $(1+J_1)$ is invertible as an operator in $\lb$ since:

\be{o(1)}
 \|J_1\|_{\lb\mapsto\lb}<\sup_{p>\epsilon}\|\hat Q_1^{-1}\|\left(\|\hat
B\|\|1\|_b+n\max_{1\le j\le n}\|{\bf D_j}\|_b\right)\rightarrow 0 
\mbox{\ as\ }b\rightarrow\infty
\ee

Given $C$, $\bfY_1$ is therefore uniquely determined from
 (\ref{eqabstract}) as an $\lb(\RR^+)$ function.

The analytic structure of $\bfY_1$ for small $z$ is contained in
in (\ref{solgen}), (\ref{solgen}'). As a result,

\be{formalseries}
\lap(\bfY_1)(x)\sim
C\sum_{k=0}^\infty\frac{\Gamma(k-\beta)}{x^{k-\beta}}\bfa_k
\ee

\z where $\sum_{k=0}^{\infty}\bfa_k z^k$ is the series of $\bfa(z)$
near $z=0$.

Correspondingly, we write (\ref{eqM})
as 

\be{eqabstractm}
(1+J_k)\bfY_k=\hat Q_k^{-1} \bfR_k
\ee

\z with $\hat Q_k:=(-\hat\Lambda+p+k)$ and 

\be{defjm}
(J_k\bfh)(p):=\hat Q_k^{-1}\left(\hat
B\int_0^p\bfh(s)ds-\sum_{j=1}^n
\int_0^ph_j(s)\bfd_j(p-s)ds\right)
\ee

\be{estiunif}
 \|J_k\|_{\lb\mapsto\lb}<\sup_{p\ge 0}\|\hat Q_k^{-1}\|\left(\|\hat B\|\|1\|_b+n\max_{1\le j\le n}\|{\bf D_j}\|_b\right)
\ee

\z Since $\sup_{p\ge 0}\|\hat Q_k^{-1}\|\rightarrow 0$ as $k\rightarrow\infty$
we have

\be{normj}
 \sup_{k\ge 1}\left\{\|J_k\|_{\lb\mapsto\lb}\right\}\rightarrow 0\ \mbox{as\ } b\rightarrow\infty
\ee

 Thus,

\Pp{Uniformnorm}
 
For large $b$, $(1+J_k), k\ge 1$ are simultaneously invertible in $\lb$,
(cf. \ref{normj}). For specified $\bfY_0$ and $C$, $\bfY_k, k\ge 1$ are
uniquely determined and moreover, for $k\ge 2$,

\be{estiY}
\|\bfY_k\|_b\le \frac{\sup_{p\ge 0}\|\hat Q_k^{-1}\|}
{1-\sup_{k\ge 1}\|J_k\|_{\lb\mapsto\lb}}\|\bfR_k\|_b:=K\|\bfR_k\|_b
\ee

\ePp

\Box

(Note: As we will see later, there only is a {\em one-parameter}
freedom in $\bfY_k$: a change in $\bfY_0$ can be compensated
by a corresponding change in $C$.) 

Because of condition $\sum m=k$ in the definition of $\bf R_{\rm
k}$, we get, by an easy induction, the homogeneity relation
with respect to the free constant $C$,

\be{homog}
\bfY_k^{[C]}=C^k\bfY_k^{[C=1]}=:C^k\bfY_k
\ee

\Pp{Prest}
 For any $\delta>0$ there is a large enough  $b$,
so that

\be{estimunifk}
\|\bfY_k\|_b<\delta^k,\ k=0,1,..
\ee

\z Each $\bfY_k$ is  Laplace transformable 
and $\bfy_k=\lap(\bfY_k)$ solve 
(\ref{systemform}). \ePp

{\em Proof}

We first show inductively
that the $\bfY_k$ are bounded. Choose $r$ small
enough and $b$ large so that $\|\bfY_0\|_b<r$.
 Note that in the expression
of $\bfR_k$, only $\bfY_i$ with $i<k$ appear. We show by induction that $\|\bfY_k\|_b<r$
for all $k$.
Using (\ref{estiY}), (\ref{eqM}) the explanation to
(\ref{systemform}) and Proposition~\ref{estimd} we get

\be{finesti1}\|\bfY_k\|_b<K\|\bfR_k\|_b\le
\sum_{|\bfl|> 1}\mu^{|\bfl|}
 r^k\sum_{\Sigma m=k}1\le r^k\left(\sum_{l>1}{l\choose k}\mu^l
\right)^n\le (r(1+\mu)^n)^k<r
\ee

\z if $r$ is small which completes this induction step. But now if we look
again at (\ref{finesti1}) we see that in fact $\|\bfY_k\|_b\le
(r(1+\mu)^n)^k$. Choosing $r$ small enough, (and to that end, $b$ large
enough) the first part of Proposition~\ref{Prest} follows. Laplace
transformability as well as the fact that $\bfy_k$ solve
(\ref{systemform}) follow immediately from (\ref{estimunifk}) (observe
again that, given $k$, there are only finitely many terms in the sum in
$\bfR_k$).

\Box

Therefore,

\Rm{r21}
The series 

\be{defgenlap}
\sum_{k=0}^{\infty}C^k(\bfY_k\cdot{\mathcal{H}})\circ\tau_k
\ee

\z is   convergent in
$\lb$ for large $b$ and thus
the sum is Laplace transformable. By Remark~\ref{contlap}
and Proposition~\ref{estimunifk}

\be{laptrfin}
\lap\left(\sum_{k=0}^{\infty}C^k(\bfY_k{\mathcal{H}})\circ\tau_{k}\right)(x)
=\sum_{k=0}^{\infty}e^{-kx}\lap(\bfY_k)(x)
\ee

\z is  uniformly convergent for
large $x$  (together with its derivatives 
with respect to $x$). Thus
(by its formal construction) (\ref{laptrfin}) is  a solution
of  (\ref{eqor}).

\eRm

\Box

(Alternatively, we could have checked in a straightforward way that the
series (\ref{defgenlap}), truncated to order $N$ is a solution of the
convolution equation (\ref{eqil}) on the interval $p\in[0,N)$ and in
view of the $\lb(\RR^+)$ (or even $\lloc$) convergence it has to be one
of the general solutions of the convolution equation and therefore
provide a solution to (\ref{eqor}).)

{\em Proof of Proposition~\ref{basicpp}, ii)}

We now show (\ref{analiticstructn}). This is done
from the system (\ref{eqM}) by induction on $k$.
For $k=0$ and $k=1$ the result follows from
Proposition \ref{ansector} and 
Proposition \ref{uniformlb}. For the induction step
we consider the operator $J_k$ (\ref{defjm}) on the space

\be{funspace3}
{\cal T}_k=
\left\{{\bf Q}:[0,\epsilon)\mapsto\CC:
{\bf Q}(z)= z^{k\beta-1}\bfA_k(z)
\right\} 
\ee

\z where $\bfA_k$ extends as an analytic function
in a neighborhood $\calv$ of $z=0$. Endowed with  the norm

$$\|{\bf Q}\|_{\cal
T_{\rm k}}:=\sup_{z\in\calv}|\bfA_k(z)|_\wedge$$

\z ${\cal T}_k$ is a Banach space. 
\Rm{analytprop}

For $k\in\NN$ the operators $J_k$ in (\ref{defjm})   
extend  continuously to ${\cal T}_k$ and their norm is $O(\epsilon)$. The functions
${\bf R}_k$, $k\in\NN$ 
(cf. (\ref{eqabstractm}), (\ref{eqM})), belong to ${\cal T}_k$.
Thus for $k\in\NN$, $\bfY_k\in {\cal T}_k$.
\eRm

\z If $A,B$ are analytic then for $z<\epsilon$

\be{ident12}
\int_0^z ds\,s^{k\beta-1}A(s)B(z-s)=z^{k\beta}\int_0^1 dt\,t^r
A(zt)B(z(1-t))
\ee

\z is in $\cal T_{\rm k}$ with norm  $O(\epsilon)$ and the assertion
about
$J_k$ follows easily. Therefore $\bfY_k\in\cal T_{\rm k}$ 
if $\bfR_k\in\cal T_{\rm k}$. We prove both these properties by induction
and (by the  homogeneity of $\bfR_k$ and the fact that
$\bfR_k$ depends only on $\bfY_m,m<k$)this  amounts 
 to checking that
if ${\bf Y}_m\in {\cal T}_m$ and ${\bf Y}_n\in {\cal T}_n$
then

$${\bf Y}_m*{\bf Y}_n\in{\cal T}_{m+n}$$

This follows from the identity

$$\int_0^z ds\, s^rA(s)(z-s)^qB(z-s)=z^{r+q+1}\int_0^1 dt\,t^r(1-t)^q A(zt)B(z-zt)$$

\Box

It is
now easy to see that $\lap_\phi\bor\hatby_k\sim\hatby_k$
(cf. Theorem \ref{teo2}). Indeed, note that in view of
Remark~\ref{analytprop} and Proposition~\ref{Prest},  $\lap(\bfY_k)$
have asymptotic power series that can be differentiated
for large $x$ in the positive half plane. Since
$\lap(\bfY_k)$ are true solutions of the system (\ref{systemform})
their asymptotic series are formal solutions of 
(\ref{systemform}) and by  the uniqueness of the formal
solution of (\ref{systemform}) once $C$ is given, the property
follows.

In the next subsection, we prove  that the general solution of the system 
(\ref{systemform}) can be obtained by means
of Borel transform of formal series and analytic continuation.

We define $\bfY^+$ to be the
function
defined in Proposition~\ref{uniright},
extended in $\cald\cap\CC^+$
 by the unique solution
of (\ref{eqil})  $\bfY_0$ provided by Proposition~\ref{ansector}.
(We define $\bfY^-$ correspondingly.)

By Proposition~\ref{uniright}, $ii)$ $\bfY^\pm$ are solutions
of (\ref{eqil}) on $[0,\infty)$
(cf. (\ref{defT})). By Lemma \ref{l0,1+eps} any solution
on $[0,\infty)$
can be obtained from, say, $\bfY^+$ by choosing $C$
and then solving uniquely (\ref{eqi3}) on $[1+\epsilon,\infty)$
(Proposition~\ref{uniright}). We now show that
the solutions of (\ref{eqabstract}), (\ref{eqabstractm})
are continuous boundary values  of functions
analytic in a region bounded by $\RR^+$.

\Rm{r10} The function $\bfd(s)$
defined in (\ref{defderiv}) by substituting $\bfH=\bfY^{\pm}$, is in
${\cal T}^{\pm}_{0,\infty}$ (cf. (\ref{defT})).

\eRm

By Proposition \ref{uniright}, $ii)$ it is easy to check that if ${\bf
H}$ is any function in $ {\cal T}^+_{0,A}$ then $\bfY^+*{\bf Q}\in {\cal
T}^+_{0,A}$. Thus, with $\bf H=Y^+$, ll the terms in the infinite sum
in (\ref{defderiv}) are in ${\cal T}^+_{0,A}$.  For fixed $A>0$, taking
$b$ large enough, the norm $\rho_b$ of $\bfY^+$ in $\lb$ can be made
arbitrarily small uniformly in all rays in $S^+_{0,A}$ (\ref{defT})
(Proposition \ref{uniright}). Then by Corollary~\ref{corolla} and
Proposition \ref{uniright} $ii)$, the uniform norm of each term in the
series (\ref{defderiv}) can be estimated by
$Const\,\rho_b^{|\bfl-1|}\nu^{|\bfl|} $ and thus the series converges
uniformly in ${\cal T^+}_{0,\infty}$, for large $b$.

\Box

\Lm{lallsol}
i) The system (\ref{eqM}) with $\bfY_0=\bfY^+$ (or $\bfY^-$)
 and given $C$ (say $C=1$)
has a unique solution  in $\lloc(\RR^+)$, namely
$\bfY_k^+$, ($\bfY_k^-$, resp.), $k\in\NN$. Furthermore,
for large $b$ and all $k$, $\bfY_k^+\in{\cal T}^+_{0,\infty}$
($\bfY_k^-\in{\cal T}^-_{0,\infty}$) (cf. (\ref{defT})).

ii) The general solution of the equation
(\ref{eqil}) in $\lloc(\RR^+)$ can be written in either
of the forms:

\ba{transssum}
&&\bfY^+ +\sum_{k=1}^{\infty}C^k(\bfY_k^+\cdot{\cal H})\circ\tau_k \cr
&&\bfY^- +\sum_{k=1}^{\infty}C^k(\bfY_k^-\cdot{\cal H})\circ\tau_k
\end{eqnarray}

\eLm

{\em Proof.}

i) The first part follows from the same arguments as
Proposition~\ref{Uniformnorm}. For the last statement
 it is easy to see (cf. (\ref{ident12})) that 
$J_k{\cal T}^+_{0,\infty}\subset {\cal T}^+_{0,\infty}$
and by Proposition~\ref{comutconv} the inequalities
(\ref{estiunif}), (\ref{normj}) hold for 
$\|\|_{{\cal T_{\rm 0,A}}\mapsto{\cal T_{\rm 0,A}}}$
($A$ arbitrary)
replacing $\|\|_{\lb\mapsto\lb}$.

$ii)$ We already know 
that  $\bfY^+$ solves (\ref{eqM}) for $k=0$. For
$k>0$ by $i)$ $C^k\bfY_k\in{\cal T}_{0,\infty}$
 and so, by continuity, the boundary values
of $\bfY_k^+$ on $\RR^+$ solve the system (\ref{eqM})
on $\RR^+$ in $\lloc$. The rest of {\em ii)} follows
from Lemma \ref{l0,1+eps}, Proposition \ref{uniright}
and the arbitrariness of $C$ in  (\ref {transssum})
(cf. also (\ref{eqneuman}).

$\Box_{_{L_4}}$

%\phantom{asd}

\begin{subsection}{Analytic structure and averaging}
\label{analyt}

Having the general structure of the solutions of (\ref{eqil})
given in Proposition \ref{convoutside} 
and in Lemma \ref{lallsol} we can 
obtain various analytic identities. The function $\bfY_0^\pm:=\bfY^\pm$
has been  defined in the previous section.

\Pp{resurgenrel} For $m\ge 0$,

\be{resufin}
\bfY_m^-=\bfY_m^++\sum_{k=1}^{\infty}{{m+k}\choose{m}}S_\beta^k(\bfY_{m+k}^+\cdot{\cal
H})\circ\tau_k 
\ee

\ePp

\emph{Proof}.

$\bfY_0^-(p)$
is a particular solution of (\ref{eqil}). It follows from Lemma \ref{lallsol}
that the following identity holds on $\RR^+$:

\be{ident1}
\bfY_0^-=\bfY_0^++\sum_{k=1}^{\infty}S_\beta^k({\bfY_k^+}\cdot{{\cal H}})\circ\tau_k
\ee

\z since, by (\ref{struct}) and (\ref{caracty0p=1}), (\ref{ident1}) holds for
$p\in(0,2)$.

 By
Lemma~\ref{lallsol} for any $C_+$ there is a $C_-$ such that

\ba{generalresurgence}
\bfY_0^+ +\sum_{k=1}^{\infty}C_+^k(\bfY_k^+\cdot{\cal H})\circ\tau_k 
=\bfY_0^- +\sum_{k=1}^{\infty}C_-^k(\bfY_k^-\cdot{\cal H})\circ\tau_k
\end{eqnarray}

\z To find the relation $C_+$ and $C_-$ we take
 $p\in(1,2)$; we get, comparing with (\ref{ident1}):

\be{stoketran}
\bfY_0^+(p)+C_+\bfY_1(p-1)=\bfY_0^-(p) +C_-\bfY_1(p-1)\Rightarrow
C_+=C_-+S_\beta
\ee

\z whence, for any $C\in\CC$,

\be{resurgencegen}
\bfY_0^+ +\sum_{k=1}^{\infty}(C+S_\beta)^k(\bfY_k^+\cdot{\cal H})\circ\tau_k 
=\bfY_0^- +\sum_{k=1}^{\infty}C^k(\bfY_k^-\cdot{\cal H})\circ\tau_k
\ee

\z Differentiating $m$ times  w.r. to $C$ and taking $C=0$ we get

$$\sum_{k=m}^{\infty}\frac{k!}{(k-m)!}S_\beta^{k-m}(\bfY_k^+\cdot{\cal H})\circ\tau_k =m!(\bfY_m^-\cdot{\cal H})\circ\tau_m$$

\z from which we obtain (\ref{resufin}) by rearranging the terms and applying $\tau_{-m}$.

\Box

\Pp{propR1}
The functions $\bfY_k$, $k\ge 0$, are analytic in $\mathcal{R}_1$.
\ePp
\emph{Proof.}

Starting with (\ref{ident1}), 
if we take $p\in(1,2)$ and obtain:

\be{ident2}
\bfY_0^-(p) =\bfY_0^+(p)+S_\beta\bfY_1(p-1)
\ee

By Proposition~\ref{uniright} and Lemma~\ref{lallsol}
the l.h.s of (\ref{ident2}) is analytic in a
lower half plane neighborhood of $(\varepsilon,1-\varepsilon),\
(\forall\varepsilon\in(0,1))$ 
and continuous in the closure of such a neighborhood.
The r.h.s. is analytic in an
upper half plane neighborhood of $(\varepsilon,1-\varepsilon),\
(\forall\varepsilon\in(0,1))$ 
and continuous in the closure of such a neighborhood.
Thus, $ \bfY_0^-(p)$ can be analytically continued 
along a path crossing the interval $(1,2)$ from
below, i.e., $\bfY_0^{-+}$ exists and is analytic. 

Now,  in (\ref{ident1}), let  $p\in(2,3)$:
\ba{ident3}
&&
S_\beta^2\bfY_2(p-2)=\bfY_0(p)^{-}-\bfY(p)^{+}-S_\beta\bfY_1(p-1)^{+}=\cr
&&\bfY_0(p)^{-}-\bfY_0(p)^{+}-\bfY_0(p)^{-+}+\bfY_0(p)^{+}=\bfY_0(p)^{-}-\bfY_0(p)^{-+}
\end{eqnarray}

\z and, in general, taking $p\in (k,k+1)$ we get

\be{identn}
S_\beta^k\bfY_k(p-k)=\bfY_0(p)^{-}-\bfY_0(p)^{-^{k-1}+}
\ee

\z Using (\ref{identn}) inductively, the same 
arguments that we used for $p\in(0,1)$ show that
$\bfY_0^{-^{k}}(p)$
can be continued analytically in the upper half plane. Thus, we have

\Rm{any0}
The function $\bfY_0$ is analytic in $\mathcal{R}_1$. In fact,
for $p\in(j,j+1),$ $k\in\NN$,

\be{formacY}\bfY_0^{-^j+}(p)=\bfY_0^+(p)+\sum_{k=1}^jS_\beta^k\bfY_k^+(p-k)\heav(p-k)
\ee
\eRm

The relation (\ref{formacY}) follows from (\ref{identn}) and (\ref{ident1}).

$\Box_{R_{\ref{any0}}}$

Note: Unlike (\ref{ident1}), in (\ref{formacY}) the sum contains
a finite number of terms. For instance we have:

\be{example0}\bfY_0^{-+}(p)=\bfY_0^+(p)+\heav(p-1)\bfY_1^+(p-1).\ \ (\forall p\in\RR^+)\ee

%%%%

The analyticity of $\bfY_m$, $m\ge 1$ is shown inductively
on $m$,
using (\ref{resufin}) and following exactly the same course 
of proof as for $k=0$.

\Box

\Rm{analyticY} If $S_\beta=0$ then $\bfY_k$ are analytic in
$\mathcal{W}_1\cup\NN$.
\eRm

Indeed, this follows from (\ref{ident1}) 
(\ref{resufin}) and Lemma~\ref{lallsol}, $i)$

\Box

On the other hand, if $S_\beta\ne 0$,  then
all $\bfY_k$ are analytic continuations of the Borel transform of $\bfy_0$
(cf. (\ref{ident3})). This is an instance of the so-called resurgence.

\z Moreover, we can now calculate $\bfY_0^{ba}$.
By definition, (see the discussion
before Remark~\ref{Rright}) on the interval $(0,2)$,

\be{med1}
\bfY_0^{ba}=\frac{1}{2}(\bfY_0^++\bfY_0^-)=
\bfY_0^++\frac{1}{2}S_\beta(\bfY_1\,\heav)\circ\tau_{1}
\ee

\z Now we are looking for a solution of (\ref{eqil})
which satisfies the condition (\ref{med1}). By 
comparing with Lemma~\ref{lallsol}, which
gives the general form of the solutions of (\ref{eqil}), we get,
now on the whole positive axis,

\be{med2}
\bfY_0^{ba}=\bfY_0^++\sum_{k=1}^{\infty}\frac{1}{2^k}S_\beta^k 
(\bfY_k^+ {\cal H})\circ\tau_k\  (\mbox{on}\ \RR^+)
\ee

\z which we can rewrite using (\ref{identn}):

\be{med3}
\bfY_0^{ba}=\bfY_0^++\sum_{k=1}^{\infty}\frac{1}{2^k}\left(\bfY_0^{-^{k}}-\bfY_0^{-^{k-1},+} \right)({\cal H}\circ\tau_k)
\ee

\Pp{pmed} 
Let $y_1(p),y_2(p)$ be analytic in $\mathcal{R}_1$,
and such  that for any path $\gamma=t\mapsto t\exp(i\phi(t))$
in $\mathcal{R}_1$, 

\be{condl1}
|y_{1,2}(\gamma(t))|<f_\gamma(t)\in\lloc(\RR^+)
\ee

Assume further that
for some large
enough $b,M$ and any path $\gamma$ in $\mathcal{R}_1$:

\be{condexpinfi}
\int_\gamma |y_{1,2}|(s)e^{-b|s|}|ds|<M
\ee

\z Then the analytic continuation $AC_\gamma(y_1*y_2) $ 
along a path $\gamma$
in $\mathcal{R}_1$,
of their
convolution product $y_1*y_2$ (defined for small $p$ by
(\ref{defconv})) exists, is locally integrable 
and satisfies (\ref{condl1}) and,
 for the same $b$  and some $\gamma$-independent
$M'>0$, 

\be{condexpinfi2}
\int_\gamma |y_1*y_2|(s)e^{-b|s|}|ds|<M'
\ee

\ePp

\emph{Proof.}

\z Since 

\be{greekmult}
2y_1*y_2=(y_1+y_2)*(y_1+y_2)-y_1*y_1-y_2*y_2
\ee

\z it is enough to take $y_1=y_2=y$. 
For $p\in\RR^+\backslash\NN$ we write:

\be{decom1}
y^-=y^++\sum_{k=1}^{\infty}(\heav\cdot y_{k}^+)\circ\tau_k
\ee

\z 
 The functions
$y_k$ are \emph{defined} inductively (the superscripts
``+,(-)'' mean, as before, the analytic
continuations in $\mathcal{R}_1$ going
below(above) the real axis). In the same way (\ref{identn})
was obtained we get by induction:

\be{ident4}
y_k=(y^{-}-y^{-^{k-1}+})\circ\tau_{-k}
\ee

\z where the equality  holds on $\RR^+\backslash\NN$
and $+,-$ mean the upper and lower continuations. 
For any $p$ only finitely many terms
 in the sum in (\ref{decom1}) are nonzero. 
The sum is also convergent in $\|\|_b$ (by dominated convergence;
note that, by assumption,
the functions $y^{--..-\pm}$ belong to the same
$\lb$).

If $t\mapsto\gamma(t)$ in $\mathcal{R}_1$, is a straight line, other
than $\RR^+$, then:

\be{forac}
AC_\gamma((y*y))=AC_\gamma(y)*_\gamma AC_\gamma(y)
 \ \mbox{if $\arg(\gamma(t))$=const$\ne0$}
\ee

\z  (Since
$y$ is analytic along such a line). The notation $*_\gamma$
means (\ref{defconv}) with $p=\gamma(t)$.

Note though that,
suggestive as it might be, (\ref{forac}) is \emph{incorrect}
if the condition stated there is not satisfied 
and $\gamma$ is a path that crosses the real line\,(see the Appendix, Section~\ref{noncomut})!

We get from (\ref{forac}), (\ref{decom1}) 
(see also (\ref{u4}), in the Appendix):

\ba{rel12}
&&(y*y)^-=y^-*y^-=
y^+*y^+ +\sum_{k=1}^{\infty}\left({\cal H}\sum_{m=0}^k
y_{m}^+*y_{k-m}^+\right)\circ\tau_k=\cr
&&(y*y)^++\sum_{k=1}^{\infty}\left({\cal H}\sum_{m=0}^k
\big(y_{m}*y_{k-m}\big)^+\right)\circ\tau_k
\end{eqnarray}

\z and now the analyticity of $y*y$ in $\mathcal{R}_1$ follows:
on the interval $p\in(m,m+1)$ we have from (\ref{ident4}) 

\be{anc1}
\big(y*y\big)^{-^j}(p)=\big(y*y\big)^-(p)=\big(y^{*2}\big)^+(p) +\sum_{k=1}^{j}\sum_{m=0}^k
\big(y_{m}*y_{k-m}\big)^+(p-k)
\ee

\z Again, formula (\ref{anc1}) is useful for analytically continuing $(y*y\big)^{-^j}$
along a path as the one depicted in Fig.1.
By dominated convergence, $(y*y)^{\pm}\in\mathcal{T}_{(0,\infty)}^{\pm}$,
(\ref{defT}).
By (\ref{ident4}), $ y_{m}$ are analytic in $\mathcal{R}_1^+:=
\mathcal{R}_1\cap\{p:\Im(p)>0\}$ and thus by (\ref{forac})
 the r.h.s. of (\ref{anc1})
can be continued analytically in $\mathcal{R}_1^+$. The  same is
then  true for  $(y*y)^{-}$. The function
$(y*y)$ can be extended analytically along paths that 
cross the real line from below. Likewise,
$(y*y)^{+}$ can be continued analytically in the lower half plane so that
$(y*y)$ is analytic in $\mathcal{R}_1$. 

Combining (\ref{anc1}), (\ref{forac}) and (\ref{greekmult}) 
we get a similar formula for the analytic
continuation of the convolution product of two functions, $f,g$
satisfying the assumptions of Proposition~\ref{pmed}

\be{anc12}
(f*g)^{-^j+}=f^+*g^+ +\sum_{k=1}^{j}\left(\heav\sum_{m=0}^k
f_{m}^+*g_{k-m}^+\right)\circ\tau_k
\ee

Note that (\ref{anc12}) corresponds to (\ref{decom1}) and
in those notations
we have:

\be{componf*g}
\big(f*g\big)_k=\sum_{m=0}^k
f_{m}*g_{k-m}
\ee

Integrability as well as (\ref{condexpinfi2})  follow
from (\ref{ident4}), (\ref{anc1}) and Remark~\ref{1}.

$\Box_{P\ref{pmed}}$

By (\ref{defmed}) and (\ref{ident4}), 

$$y^{ba}=y^++\sum_{k=1}^\infty\frac{1}{2^k} (y_{k}^+\heav)\circ\tau_k$$

\z so that (see (\ref{u4}))

\ba{commutproof}&&y^{ba}*y^{ba}=\Big(y^++\sum_{k=1}^{\infty}\frac{1}{2^k}
({\cal H}\circ\tau_k)( y_{k}^+\circ\tau_k)\Big)^{*2}=\cr
&&y^+*y^++\sum_{k=1}^{\infty}\frac{1}{2^k}{\cal H}\circ\tau_k\sum_{m=0}^k
(y_{m}^+\circ\tau_m)*(y_{k-m}^+\circ\tau_{k-m})\circ\tau_k=\cr
&&y^+*y^++\sum_{k=1}^{\infty}\frac{1}{2^k}{\cal H}\circ\tau_k\sum_{m=0}^k
(y_{m}*y_{k-m})^+\circ\tau_k=(y^{*2})^{ba}
\end{eqnarray}

To finish the proof of Theorem \ref{teo2}
note that on any finite interval the sum in (\ref{defmed}) has
only a  finite number of terms
and by (\ref{commutproof})  balanced averaging  commutes with any finite sum of the type
\be{finisum}
\sum_{k_1,..,k_n}c_{k_1..k_n}f_{k_1}*..*f_{k_n}
\ee

\z and then, by continuity, with any sum
of the form (\ref{finisum}), with a finite or infinite number of terms, 
provided it converges in $\lloc$. Averaging
thus commutes with all the operations
involved in the equations (\ref{eqabstractm}).
By uniqueness therefore, if $\bfY_0=\bfY^{ba}$ 
then  $\bfY_k=\bfY_k^{ba}$ for all $k$. Preservation of reality is
immediate since (\ref{eqil}), (\ref{eqM}) are real if (\ref{eqor}) is real, therefore
$\bfY_0^{ba}$ is real-valued on $\RR^+\backslash\NN$ (since it is real-valued on
$[0,1)\cup(1,2)$) and so are, inductively, all $\bfY_k$.

\end{subsection}

\end{subsection}
\begin{section}{Acknowledgments}

The author would like
to thank Professors Michael Berry, Percy Deift and  Jean Ecalle
for very interesting discussions. Special thanks are due to 
Prof. Martin Kruskal for long and illuminating debates
and, together with Prof. Joel Lebowitz, for their
continued support and encouragements throughout this work.
\end{section}

\setcounter{section}{0}
\renewcommand{\thesection}{\Alph{section}}

\begin{section}{Appendix}
\begin{subsection}{Example of  non-typical behavior}

Consider the equation

\be{exa1}
f'=-f-\frac{1}{2x}f+
\frac{1}{x}-\frac{1}{2x^2}
\ee

\z The general solution of this equation is given by

\be{ape2}
f=\frac{1}{x}+Cx^{-1/2}e^{-x}=\int_0^{\infty}\left(p+\frac{C}{\sqrt{p-1}}{\cal H}(1-p)\right)dp
\ee

\z We see that the asymptotic series of $f$
for $x\rightarrow \infty,\ \Re(x)>0$, $\hatby_0=1/x$.
 The inverse Laplace transform of $f$ is 

\be{ape3}
\lapi f=p+\frac{C}{\sqrt{p-1}}{\cal H}(1-p)
\ee

i)
 The Stokes constant is zero and $\bfY_0=\bor(\hatby_0)=p$ is entire.

ii) All combinations $\lambda\bfY_0^++(1-\lambda)\bfY_0^-$
coincide. Therefore (\ref{combilin}) does not hold.

 Equation (\ref{exa1})  is exceptional, in the sense that the properties $i), ii)$
above do not withstand 
a small perturbation. Indeed, for the equation

\be{eqpert}f'=-f-\frac{1}{2x}f+
\frac{1+\epsilon}{x}-\frac{1}{2x^2}
\ee

\z we have  $\bor(\hatby_0)=2\epsilon+p+\epsilon(1-p)^{-1/2}$
and the inverse Laplace transform of the general solution is

$$\lapi(f)=\left\{ \begin{array}{cc}
2\epsilon+p+\epsilon(1-p)^{-1/2}\ &\mbox{for $p<1$}\cr
2\epsilon+p+C(p-1)^{-1/2}\ &\mbox{for $p>1$}\cr
\end{array}\right.
$$

\end{subsection}
\begin{subsection}{$AC(f*g)$ versus $AC(f)*AC(g)$}
\label{noncomut}

Typically, the  analytic continuation along curve in ${\cal W}_1$
which is not homotopic to a straight line will not commute
with convolution.
 For example, in
equation (\ref{eqpert}),  $\bor(\hatby_0)^{-+}* \bor(\hatby_0)^{-+}
\ne[\bor(\hatby_0)*\bor(\hatby_0)]^{-+}$, as it 
can be  seen from Remark~\ref{rapen}
below (or by direct calculation). This situation
is generic:

\Rm{rapen}
Let $y$ be a function satisfying the conditions
stated in Proposition~\ref{pmed} and assume that
$p=1$ is a branch point of $y$. Then,

\be{contrex2}
(y*y)^{-+}\ne y^{-+}*y^{-+}
\ee
\eRm

{\em{Proof}}

\z Indeed,
by (\ref{anc12}) and  (\ref{ident4})

\ba{contrexgen}
&&(y*y)^{-+}=y^+*y^++2[(y^+*y_1^+)\heav]\circ\tau_1
\ne y^{-+}*y^{-+}=\cr&&
[y^++(\heav y_1^+)\circ\tau_1]^{*2}=
y^+*y^++2[(y^+*y_1^+)\heav]
\circ\tau_1+[\heav(y_1^+*y_1^+)]\circ\tau_2\cr&&
\end{eqnarray}
\z since in view of   (\ref{ident4}), in our assumptions,
 $y_1\not\equiv 0$ and thus $y_1*y_1\not\equiv 0$.

\Box

There is also the following intuitive reasoning leading to the same
conclusion. For a generic system of the form
(\ref{eqor})--(\ref{asymbh}), $p=1$ is a branch point of $\bfY_0$ and so
$\bfY_0^-\ne\bfY_0^{-+}$. On the other hand, if $AC_{-+}$ commuted with
convolution, then $\lap(\bfY_0^{-+})$ would provide a solution of
(\ref{eqor}).  By Lemma~\ref{lallsol}, $\lap(\bfY_0^-)$ is a different
solution (since $\bfY_0^-\ne\bfY_0^{-+}$). As $\bfY_0^-$ and
$\bfY_0^{-+}$ coincide up to $p=2$ we have
$\lap(\bfY_0^{-+})-\lap(\bfY_0^-)=O(e^{-2x}x^{power})$ for $x\rightarrow
+\infty$.  By Theorem~\ref{teo2} however, no two solutions of
(\ref{eqor})--(\ref{asymbh}) can differ by less than $e^{-x}x^{power}$
without actually being equal (also, heuristically, this can be checked
using formal perturbation theory), contradiction.

\end{subsection}
\begin{subsection}{Useful formulas}
\label{usefulfor}

\be{u1}{\cal B}(\frac{1}{x^n})=\frac{p^{n-1}}{\Gamma(n)}
\ \mbox{or}\ {\cal L}(p^{n})=\frac{\Gamma(n+1)}{x^{n+1}}\ee

\be{u3}p^q*p^r=\frac{\Gamma(q+1)\Gamma(r+1)}{\Gamma(q+r+2)}
p^{q+r+1}\ee
\ with $f_{1,2}(p):=p\mapsto\heav(p-k_{1,2})g_{1,2}(p-k_{1,2})$
we have
\be{u4}\Big(f_1*f_2\Big)(p)=\heav(p-k_1-k_2)\Big(g_1*g_2\Big)(p-k_1-k_2)\ee

\end{subsection}
\end{section}

\end{document}